\documentclass[12pt]{article}
\usepackage{epsfig}
\usepackage{amsmath}
\usepackage{amssymb}

\newcommand{\no}{\nonumber}
\newcommand{\be}{\begin{equation}}
\newcommand{\ee}{\end{equation}}
\newcommand{\bi}{\begin{itemize}}
\newcommand{\ei}{\end{itemize}}
\newcommand{\br}{\begin{eqnarray}}
\newcommand{\er}{\end{eqnarray}}

\newcommand{\qed}{\mbox{$\square$}\newline}

\newcommand{\abs}[1]{\lvert #1 \rvert}
\newcommand{\norm}[1]{\lVert #1 \rVert}

\newtheorem{theo}{Theorem}[section]
\newtheorem{defin}{Definition}[section]
\newtheorem{prop}{Proposition}[section]
\newtheorem{lem}{Lemma}[section]
\newtheorem{cor}{Corollary}[section]

\begin{document}

\title{Variational Principle of KPP Front Speeds \\in Temporally Random Shear Flows \\
with Applications}

\author{James Nolen\thanks{Department of Mathematics, University of Texas at Austin,
Austin, TX 78712 (jnolen@math.utexas.edu).} \;
and \; Jack Xin\thanks{Department of Mathematics, University of California at Irvine, Irvine, CA 92697.
(jxin@math.uci.edu).}}

\date{}
%\date{\today}
\maketitle
\begin{abstract}
We establish the variational principle of Kolmogorov-Petrovsky-Piskunov (KPP) front speeds
in temporally random shear flows inside an infinite cylinder, under suitable assumptions
of the shear field. A key quantity in the variational principle is the
almost sure Lyapunov exponent of a heat operator with random potential.
The variational principle then allows us to bound and compute the
front speeds. We show the linear and quadratic laws of speed enhancement as well as a
resonance-like dependence of front speed on the temporal shear correlation length. To
prove the variational principle, we use the comparison principle of solutions,
the path integral representation of solutions,
and large deviation estimates of the associated stochastic flows.

\end{abstract}

\date{}
\setcounter{page}{1}
\section{Introduction}
\setcounter{equation}{0}
Reaction-diffusion front propagation in strongly time dependent random media arises
in premixed flame propagation problems (\cite{CW,MS2,Peters,Ro,Xin3,Yak} and references),
interacting particle systems (\cite{MS,CD} and references) and population biology
(\cite{Shen} and references). A fundamental issue is
to characterize, bound and compute the large time front speed, an upscaled
quantity that depends on statistics of the random media in a highly nontrivial manner.
In combustion literature, ad hoc and formal procedures abound for
approximation, such as closures and renormalization group methods \cite{Peters,Yak}.
In this paper, we establish variational principles of propagation speeds of
KPP reaction-diffusion fronts through
temporally random shear flows inside an infinite cylinder.
The variational characterization then allows us to estimate and compute the statistical
properties of front speeds with both accuracy and ease.

The model equation is:
\be
u_t =\frac{1}{2} \Delta_{z} u + B \cdot \nabla_{z} u + f(u), \label{e0}
\ee
where $u=u(z,t)$, $z=(x,y) \in R \times \Omega$, $\Omega$ a bounded open
subset of $R^{n-1}$ with Lipschitz continuous boundary, $n \geq 2$;
$B = (b(y,t), 0, \dots, 0)$, $b(y,t)$ is a stationary Gaussian process in $t$, with a given
deterministic profile in $y$, to be made more precise later.
The nonlinear function $f(u) \in C^{1}([0,1])$ is a KPP nonlinearity:
$f(u) >0$ for $u \in (0,1)$, $f(0)=f(1)=0$,
$f'(0) = \sup_{u \in (0,1)} f(u)/u$. An example is $f(u)=u(1-u)$.

For compactly supported initial data bounded between $0$ and $1$, solutions of (\ref{e0})
develop into propagating fronts separating the cylindrical domain
into a region where $u \approx 1$ and the rest where $u \approx 0$, which
correspond to burned (hot) and unburned (cold) states in combustion.
In case $B$ is periodic in $z$ and $t$, KPP type front dynamics and speeds
have been recently studied
for both shear and more general incompressible flows \cite{Kh,MS1,MS2,NX1,nolen:ekt04,NRX}.
Exact traveling front solutions exist \cite{NX1,nolen:ekt04,NRX},
extending those in spatially periodic media,
\cite{BH1, berestycki:tfc92, xin:epf92, xin:est91},
see also \cite{B2} and \cite{Xin1} for reviews.

For temporally random shear flows, it is more efficient to study
front solutions asymptotically without constructing exact traveling fronts.
This line of work goes back to Freidlin \cite{FR1}
where variational principles of KPP front speeds in spatially periodic
media are obtained by combining the large
deviation techniques and Feynman-Kac representation formulas of KPP solutions.
We shall further develop this approach to treat the temporally random shear flows which
generate more complexities in path integrals and unbounded variations in time.

Let us make precise our assumptions on the shear field. The function
$b(y,t) = b(y,t,\hat \omega)$ is a mean zero Gaussian random field over $(y,t)$,
periodic in $y$ with period $L$ for each fixed $t$, and
stationary in $t$ for each fixed $y$.
The field $b$ is defined over probability space $(\hat \Omega, \hat{\mathcal{F}}, Q)$
and has covariance function $\Gamma(y_1,y_2, t_1,t_2) = E_Q[b(y_1,t_t)b(y_2,t_2)]$.
The following assumptions hold on $b(y,t)$:
\bi
\item[{\bf A1:}] (Periodicity in $y$) Let $C^{0,1}_P(D)$ denote the space of Lipschitz continuous functions
that are periodic on the period cell $D = [0,L]^{n-1}$. For each $\hat \omega \in \hat \Omega$,
there is a continuous map $J( \cdot, \hat \omega): [0,+\infty) \to C^{0,1}_P(D)$ such that
$b(\cdot , t , \hat \omega) = J(t, \hat \omega)$.
\item[{\bf A2:}] (Stationarity in $t$) For each $s \in R^+$ there is a measure preserving transformation $\tau_s:\hat \Omega \to \hat \Omega$ such that $b(y, \cdot + s, \hat \omega) = b(y,\cdot,\tau_s \hat \omega)$. Hence, $\Gamma$ depends only on $y_1,y_2$ and $\abs{t-s}$.
\item[{\bf A3:}] (Ergodicity) The transformation $\tau_s$ is ergodic: if a set $A \in \hat {\mathcal{F}}$ is invariant under the transformation $\tau_s$, then either $Q(A) = 0$ or $Q(A) = 1$.
\item[{\bf A4:}] The field $b$ is mean zero, almost surely continuous in $(y,t)$, and has uniformly bounded variance:
\be
E[b(y,t)] = 0 \;\;\; E[b(y,t)^2] \leq \sigma^2 \;\;\; \text{for all} \;\; y \in D, t \geq 0.
\ee
\item[{\bf A5:}] The function $\hat \Gamma(r) = \sup_{y_1,y_2} \Gamma(y_1,y_2,0,r)$ is integrable over $[0,\infty)$:
\be
\int_0^\infty \hat \Gamma(r)\,dr = p_1 < \infty
\ee
for some finite constant $p_1 > 0$.
\item[{\bf A6:}] There is a finite constant $p_2 > 0$ such that
\be
\abs{\Gamma(y_1,y_2,s,t) - \Gamma(y_1,y_3,s,t) }\leq p_2 \abs{y_3 - y_2} \hat \Gamma(\abs{s-t}). \no
\ee

\ei
For example, a field satisfying Assumptions A1-A6 might have the form $b(y,t) = \sum_{j=1}^N b_1^j(y) b_2^j(t)$, where the functions $b_1^j(y)$ are deterministic, Lipschitz continuous and periodic over $D$, and the functions $b^j_2(t)$ are mean zero, stationary Gaussian fields in $t$.

Before stating the main results, let us define the family of Markov processes associated
with the linear part of the operator in (\ref{e0}). For a fixed
$\hat \omega \in \hat \Omega$ and for each $z \in R^n$, $t \geq 0$,
let $Z^{z,t}(s) = (X^{z,t}(s), Y^{z,t}(s)) \in R^n$ solve the It\^o equation:
\br
dX^{z,t}(s) &= & dW_1(s) + b(W_2(s),t-s)\,ds \no \\
dY^{z,t}(s) & = & dW_2(s) \label{e2}
\er
with initial condition $Z^{z,t}(0) = z =(x,y) \in R^n$, where
$W(s) = (W_1(s),W_2(s)) \in R^n$ is the $n$-dimensional Wiener process
with $W(0) = 0$. Let $P^{z,t}$ denote the corresponding family of measures on $C([0,t];R^n)$.
As we will see, the KPP front speed depends on large deviations of
the random variable
\be
\eta^t_z(\kappa t) = \frac{x - X^{z,t}(\kappa t) }{\kappa t}
\ee
which is the first component of the average speed of a trajectory over
time interval $[0, \kappa t]$. We consider only the first component since we are
concerned only with propagation in the $x$ direction.
The need for the parameter $\kappa$ results from the time
dependence of the field $b(y,t)$ and will become more apparent later.

Now we state the main results. First, the following lemma allows us to characterize the speed of propagation:
\begin{lem} \label{lem:muexist}
Assume that A1-A6 hold for the process $b(y,t)$. Then for each $\lambda \in R$, the limit
\be
\mu(\lambda, z) = \mu(\lambda) = f'(0) + \lim_{t \to \infty} \frac{1}{ t} \log E \left[ e^{ -\lambda X^{z,t}(t)}\right] \label{e25}
\ee
exists and is a finite constant, almost surely with respect to measure $Q$ and  independent of $z \in R^n$. Moreover, $\mu(\lambda)$ is a convex, positive, and even function of $\lambda$. Also,  $\mu(\lambda)/\abs{\lambda} \to + \infty$ as $\abs{\lambda} \to \infty$.
\end{lem}
If we let $H(c)$ be the Legendre transform of $\mu(\lambda)$
\be
H(c) =  \sup_{\lambda \in R} [c \cdot \lambda -  \mu(\lambda)], \no
\ee
then we find that the speed of propagation can be bounded above in terms of $H$.
\begin{theo}[Upper bound on front speed] \label{theo:frspeed1} Let $b(y,t,\hat \omega)$ satisfy assumptions A1 - A6. Let $u(x,y,t,\hat \omega)$ solve (\ref{e0}) with initial condition $u(x,y,0,\hat \omega) = u_0(x)$, where $u_0(x) \in [0,1]$ has compact support and is independent of $\hat \omega$. Then, for any closed set $F \subset \{ c \in R |\; H(c) > 0\}$,
\be
\lim_{t \to \infty} \sup_{y \in D} u(ct, y, t, \hat \omega) = 0 \no
\ee
uniformly in $c \in F$, for almost every $\hat \omega \in \hat \Omega$.
\end{theo}
Therefore, if we define the constant $c^* > 0$ by the variational formula
\be
c^* = \inf_{\lambda > 0} \frac{\mu(\lambda)}{ \lambda} , \label{e65}
\ee
we see from the definition of $H$ that the front spreads asymptotically with speed no greater than $c^*$ in the positive $x$ direction and with speed no greater than $c^*$ in the negative $x$ direction. Although the solution $u$ depends on $\hat \omega \in \hat \Omega$ since $B$ is a random variable over $\hat \Omega$, the function $H(c)$ and the speed $c^*$ are independent of $\hat \omega$. They are almost surely constant with respect to $\hat Q$, a consequence of the ergodicity assumption A3. We will generally suppress the dependence of $u$ on $\hat \omega$ for clarity of notation.

The constant $c^*$ is also almost surely a lower bound on the speed of propagation:

\begin{theo}[Lower bound on front speed] \label{theo:frspeed} Let $b(y,t,\hat \omega)$ satisfy assumptions A1 - A6. Let $u(x,y,t,\hat \omega)$ solve (\ref{e0}) with initial condition $u(x,y,0,\hat \omega) = u_0(x)$, where $u_0(x) \in [0,1]$ has compact support and is independent of $\hat \omega$. Then, for any compact set $K \subset \{ c \in R | \; H(c) < 0 \}$,
\be
\lim_{t \to \infty} \inf_{y \in D} u(ct,y, t, \hat \omega) = 1 \no
\ee
uniformly in $c \in K$,  for almost every $\hat \omega \in \hat \Omega$.
\end{theo}

Theorems \ref{theo:frspeed1} and \ref{theo:frspeed} extend our recent results on KPP
front speeds in temporally periodic incompressible flows \cite{nolen:ekt04,NRX} and a classical
result of Freidlin (Theorem 7.3.1, p. 494 of \cite{FR1}) where he treated the
case of spatially periodic advecting flows.
Our proofs are built on his, with additional ingredients
to handle both the time-dependence and the stochastic nature of the field $B$. For example,
in the periodic case, $\mu(\lambda)$ is the principal eigenvalue of a periodic-parabolic operator
\cite{nolen:ekt04,NRX}, and
perturbation theory \cite{Kato} implies that $\mu(\lambda)$ is differentiable in $\lambda$.
It then follows from Theorem 7.1.1 and Theorem 7.1.2 of \cite{FR1} that
the random variables $\eta^t_z( t)$ satisfy a large deviation principle with
convex rate function $S(c)$ given by (\ref{e53}).
However, if $\mu(\lambda)$ is not known to be differentiable,
the large deviation property needs a new proof. In the present case,
$\mu(\lambda)$ is not an eigenvalue of a linear operator,
so we cannot readily apply the perturbation theory \cite{Kato} to get differentiability.
Instead, we will show that a rate function exists and is convex, and that
it continues to satisfy (\ref{e53}). In fact, $\mu(\lambda)$ is related to
the almost sure principal Lyapunov exponent of a heat operator with random potential \cite{Xin3},
known as the parabolic Anderson problem (\cite{CarMol,CM} and references).  Dynamical
aspects of principal Lyapunov exponents as an extension of principal eigenvalues are recently studied in
\cite{MShen}. Regularity of $\mu(\lambda)$ is an interesting problem in itself.

The paper is organized as follows. In Section 2, we prove Theorem \ref{theo:frspeed1}.
In Section 3, we adapt Freidlin's method to prove Theorem \ref{theo:frspeed},
assuming a few technical estimates as well as a crucial lower bound stated
in Lemma \ref{lem:explowbnd}. In Sections 4 and 5 we will prove the
lower bound of Lemma \ref{lem:explowbnd} and necessary large deviation estimates.
It is in these sections that we resolve the main difficulties that
result from the additional stochastic time dependence.
We will make frequent use of the subadditive ergodic theorem
and the Borell inequality for Gaussian fields \cite{Adler2,AK, King}. In Section 6, we use the variational formula (\ref{e65}) to numerically compute the front
speed $c^*$ in an example where the temporal random process is
Ornstein-Uhlenbeck (O-U). We study the dependence of front speeds on
the covariance of the random media as well as their growth laws in small and large
advection limits. The speed grows linearly in the limit of large advection, and
quadratically in the limit of small advection.
The variational principle also yields analytical
bounds on $c^*$, demonstrating the enhancement by random shear flows and that
small temporal correlation reduces the front speed. This is analogous to the
decrease of front speeds with increasing frequency of temporally oscillating
periodic shear flows \cite{Kh,NX1,nolen:ekt04}. Furthermore, there is an optimal correlation
length that maximizes the enhancement.
Linear and quadratic speed growth laws are known for deterministic
flow patterns with channel structures (\cite{B2,BH1,BHN,Const1,HPS,KR,Vlad} and references),
also for spatially random shears inside infinite cylinders \cite{NX2,NX3} or white
in time Gaussian shears in the entire space \cite{Xin3}.
The speed growth laws known to date are not sensitive in the form of nonlinearities as long as
fronts propagate out of the initial data. In this sense, KPP plays the role of solvable model
and KPP front speeds carry universal properties.
In Section 7, we make concluding remarks.
Finally in the two Appendices, we prove Lemma \ref{lem:muexist} and several key technical bounds.

Though the arguments in our proofs rely on the periodicity of $b(y,t)$ in $y$
to provide compactness in the $y$ dimensions, they can be easily modified to solve
the same problem in an infinite cylinder with the zero Neumann boundary condition on the
sides of the cylinder. The compactness property remains, and
the process $X^{z,t}(s)$ just needs to be reflected when it hits the
boundary $R \times \partial \Omega$.
It is also not necessary for the process $b(y,t)$ to be Gaussian.
Our proofs of the bounds in Sections 4 and 5 shall
rely on the powerful Borell inequality for Gaussian process,
yet it is easy to see that if the estimates of Lemma \ref{corol:taubound} and
the Appendices hold for a given process, and if assumptions A1 - A6 hold,
the main results extend.

\section{Proof of Theorem \ref{theo:frspeed1}}
\setcounter{equation}{0}
The proof of Theorem \ref{theo:frspeed1} is based on the assumption that $f(u) \leq f'(0) u$.  This allows us to construct a super-solution to equation (\ref{e0}) as follows. Note that the integral form of (\ref{e2}) is
\br
X^{z,t}(s) &=& x + \int_0^s b(W^y_2(\tau),t-\tau) \,d\tau +  W^0_1(s) \label{e12}\\
Y^{z,t}(s) &=& W^y_2(s) \no
\er
where $W^z(s)$ denotes a Wiener process starting at $W^z(0) = z$, $P$-a.s. Using (\ref{e12}), we can express (\ref{e25}) as
\br
\mu(\lambda) &=& f'(0) +  \lim_{t \to \infty} \frac{1}{ t} \log E_{z} \left[ e^{ - \lambda \int_0^{ t} b(W^y_2(s),t-s)\,ds } e^{-\lambda W^0_1(t)}\right] \no \\
& = & f'(0) + \frac{\lambda^2}{2} + \lim_{t \to \infty} \frac{1}{ t} \log E_{y} \left[ e^{ -\lambda \int_0^{t} b(W^y_2(s),t-s)\,ds } \right] . \label{asle2}
\er
Here we have used the independence of $W_1$ and $W_2$. Note that we will sometimes use $\rho(\lambda) = \mu(\lambda) - f'(0) - \lambda^2$ to refer to the limit on the right hand side of (\ref{asle2}). Using the Feynman-Kac formula, we see that $\mu(\lambda)$ is just
\be
\lim_{t \to \infty} \frac{1}{ t} \log E_{y} \left[ e^{\int_0^{t} \lambda^2/2 + f'(0) - \lambda b(W^y_2(s),t-s)\,ds } \right] = \lim_{t \to \infty} \frac{1}{ t} \log \Phi(y,t) \no
\ee
where $\Phi = \Phi(y,t) >0$ is periodic in $y$ and solves the auxiliary initial value problem
\br
\Phi_t = \frac{1}{2} \Delta_y \Phi +  (\lambda^2/2 + f'(0) - \lambda b(y,t)) \Phi. \no \\
\Phi(y,0) \equiv 1 . \no
\er

Suppose $c > c^*$ (i.e. $H(c) > 0$). Then for $\epsilon > 0$ sufficiently small, there exists $\lambda > 0$ such that $\mu(\lambda) < (c-\epsilon) \lambda$. By Lemma \ref{lem:muexist} $\mu(\lambda)$ is finite, so that
\be
\lim_{t \to \infty} \frac{1}{t} \log \Phi(y,t) = \mu(\lambda), \no
\ee
and there is a function $R = R(y,t)$ such that $\abs{R} \to 0 $ as $t \to \infty$, uniformly in $y$, such that
\be
\Phi(y,t) = e^{\mu(\lambda)t + R(y,t)t} < e^{\lambda (c-\epsilon) t + R(y,t)t}. \no
\ee
Defining $\psi(x,y,t) = e^{-\lambda x}\Phi(y,t)$, we see that $\psi$ solves the equation
\br
\psi_t = \frac{1}{2} \Delta_{x,y} \psi + b \psi_x + f'(0) \psi \no \\
\psi(x,y,0) = e^{-\lambda x}. \no
\er
By the properties of $f$, we then have
\be
\psi_t \geq \frac{1}{2} \Delta_{x,y} \psi + b \psi_x + f(\psi) \no
\ee
when $\psi \in [0,1]$. Multiplying $\psi$ by a constant $C$ if necessary, we may assume $\psi(z,0) > u_0(z)$, and by the maximum principle we must have $u(z,t) \leq \psi(z,t)$ for all $t \geq 0$. If we define the half-space
\be
\Sigma^+_r = \left\{ (x,y) \in R^n\;| \; x > r \right\} \no
\ee
we now see that
\br
\lim_{t \to \infty} \sup_{z \in \Sigma^+_{ct}} u(x,t) &\leq& \lim_{t \to \infty} \sup_{z \in \Sigma^+_{ct}} \psi(z,t) \no \\
& = & \lim_{t \to \infty} \sup_{z \in \Sigma^+_{ct}} e^{-\lambda x}\phi(y,t) \no \\
& \leq & \lim_{t \to \infty} \sup_{y \in D} e^{-\lambda c t} e^{\lambda (c - \epsilon) t + R(y,t)t} = 0 \no
\er
since $\abs{R(y,t)} < \frac{\epsilon \lambda}{2}$ for $t$ sufficiently large. A similar argument holds for propagation in the $-x$ direction.  This completes the proof of Theorem \ref{theo:frspeed1}. \qed

%The effect of such nontrivial advection $B$ is to
%greatly enhance the speed of the traveling waves or
%of the moving interface separating regions where
%$u \approx 1$ and $u \approx 0$.  As a result,
%it is of practical interest to understand how structures
%of the flow field $B$ contribute to the enhancement of the propagation.

\section{Proof of Theorem \ref{theo:frspeed}}
Proving Theorem \ref{theo:frspeed} is more difficult, and we will need the following important estimates. The first is a lower bound analogous to Lemma 7.3.3 of \cite{FR1}:
\begin{lem} \label{lem:explowbnd}
For any compact set $K \subset \{ c \in R | \; H(c) > 0 \}$,
\be
\liminf_{t\to \infty} \frac{1}{t} \log \inf_{c \in K, y \in D} u(ct,y,t) \geq - \max_{c \in K} H(c). \label{e27}
\ee
\end{lem}
The second estimate gives a coarse bound on very large excursions of the random process $X^{z,t}$:
\begin{lem} \label{corol:taubound}
There are constants $K_1, K_2 > 0$ independent of $\kappa \in (0,1]$ such that, except on a set of $ Q$-measure zero,
\be
\sup_{ z \in R^n} P\left( \sup_{s \in [0,\kappa t]} \abs{X^{z,t}(s) - x} \geq \eta t \right) \leq K_1 e^{- K_2 \eta^2 t/\kappa} \no
\ee
for any $\kappa \in (0,1]$, $\eta > 0$, and for $t$ sufficiently large depending on $\hat \omega$, $\kappa$, and $\eta$.
\end{lem}

These lemmas represent the main technical difficulty in extending the work of \cite{FR1} to the present case with a stochastic time dependence in $b(y,t)$. For the moment, however, we delay the proof of these lemmas and show how these bounds lead to Theorem \ref{theo:frspeed}. Lemma \ref{lem:explowbnd} is proved in the next section; Lemma \ref{corol:taubound} is proved in the appendix. The proof of Theorem \ref{theo:frspeed} is based on the observation that when $u < h < 1$, the reaction rate can be bounded below. For each $u \in (0,1]$, define the reaction rate $\zeta$ by
\be
\zeta(u) = \frac{f(u)}{u} \no
\ee
and $\zeta(0) = f'(0)$. Now equation (\ref{e0}) can be written
\be
u_t = \frac{1}{2}\Delta_{z} u +  b(y,t) u_x + \zeta(u)u. \label{e3}
\ee
By the properties of $f(u)$ we see that $\zeta(u) > 0$ for $u \in [0,1)$, $\zeta(u)$ is continuous for $u \in [0,1]$, and $\zeta(0) \geq \zeta(u)$ for any $u \in [0,1]$. If $h \in (0,1)$ we define a lower bound on $\zeta$:
\be
\zeta_h = \inf_{u \in (0,h)} \zeta(u) > 0. \no
\ee
So, in regions where $u$ is bounded away from one, the reaction rate can be bounded below by $\zeta_h >0$.

For a fixed $\hat \omega \in \hat \Omega$, we can estimate $u(z,t)$ using the Feynman-Kac formula for the solution of (\ref{e3}):
\be
u(z,t) = E \left [e^{\int_0^t \zeta(t-s,u(Z^{z,t}(s),t-s))ds}u_0(Z^{z,t}(t)) \right], \label{fkeq}
\ee
where the expectation is with respect to measure $P^{z,t}$. If $\tau$ is any stopping time, we also have
\be
u(z,t) = E \left [e^{\int_0^{t \wedge \tau} \zeta(t-s,u(Z^{z,t}(s),t-s))ds}u(Z^{z,t}(t-(t\wedge \tau)),t-(t\wedge \tau) ) \right], \label{fkeq2}
\ee
where $t \wedge \tau = \min(t,\tau)$. Therefore, we can obtain estimates on $u$ by carefully choosing stopping times and restricting the expectation to certain sets of paths.  The exponential term inside the expectation will be large when the path $Z^{z,t}(s)$ passes through regions where $u$ is small and the reaction rate is large. On the other hand, if $u(Z^{z,t}(t-(t\wedge \tau)),t-(t\wedge \tau) )$ is too small, then the expectation as a whole may be small.

Now we present the argument of \cite{FR1} (see p. 494).  For $s \in R$, define the set
\be
\Psi(s) = \{ c \in R | \; H(c) = s \} \;\;\text{and} \;\;\; \underline\Psi(s) = \{ c \in R | \; H(c) \leq s\} \no
\ee
and then for any $\delta > 0$ and $T > 1$,
\be
\Gamma_T = \left( [\{ 1\} \times \underline\Psi(\delta)] \cup [ \bigcup_{1 \leq t \leq T} \{ t \} \times (t \Psi(\delta))]\right) \times R^{n-1} . \no
\ee
This defines the boundary of a region that spreads outward in $x$, linearly in $t$. Outside this region $u$ may be close to zero, but on the boundary of this region, we have the crucial lower bound from Lemma \ref{lem:explowbnd}:
\be
u(x,y,s) \geq e^{-2 \delta t} \;\;\;\text{for all}\;\; (x,y,s) \in \Gamma_t \label{e24}
\ee
for $t$ sufficiently large.

For $h \in (0,1)$, and $z$ fixed, $t, \eta > 0$, define the Markov times
\br
\sigma_h(t)  & = & \min \{ s \in [0,t] | \; u(Z^{z,t}(s),t-s) \geq h \} \no \\
\sigma_\Gamma(t) & = &\min \{s \in [0,t]|\; (Z^{z,t}(s),t-s) \in \Gamma_t \} \no \\
\tau_\eta (t) &=& \min \{ s \in [0,t] | \; \abs{X^{z,t}(s) - x} > \eta t \} \no
\er
(We set these variables equal to $+ \infty$ if the sets on the right are empty.) Then $P^{z,t}(\sigma_h(t) \leq t)$ is the probability that a particle starting at $z$ will encounter the ``hot region", $u \geq h$, at or before time t. Since $\zeta \geq 0$, it is clear from (\ref{fkeq2}) that
\be
u(ct,y,t) \geq h P^{ct,y,t}(\sigma_h(t) \leq t). \label{e26}
\ee
For simplicity, we will write $P^{ct,t}$ or just $P^{ct}$ to denote $P^{(ct,y),t}$ when the $y$ and $t$ dependence is clear.
If we choose $c$ too large, we expect that this probability will be small, since $x = ct$ would be far beyond the spreading front, beyond the region where $u \geq h$. On the other hand, we want to show that the probability is large for $c \in K$ (i.e. $\abs{c} < c^*$), for if we can show that for each $h$, $P^{ct,t}(\sigma_h(t) \leq t) \to 1$ uniformly over $c \in K$, then (\ref{e26}) implies the desired result: $u(ct,y,t) \to 1$ uniformly over $c \in K$. Note that $P^{ct,t}(\sigma_h(t) > t) = P^{ct,t}(\sigma_h(t) > t, \; \tau_\eta(t) > t) + P^{ct,t}(\sigma_h(t) > t,\; \tau_\eta(t) < t)$.  By Lemma \ref{corol:taubound},
\be
 \sup_{z \in R^n} P^{z,t} \left( \tau_\eta (t) < t \right) \to 0 \label{e16}
\ee
as $t \to \infty$ except on a set of $Q-$measure zero. So, it suffices to show that
\be
P^{ct,t}(\sigma_h(t) > t, \; \tau_\eta(t) > t) \to 0
\ee
uniformly in $c \in K$. Note that in \cite{FR1}, estimate (\ref{e16}) followed from the uniform boundedness of the field $B$, a property that we do not have in the stochastic case. The rest now follows exactly as in \cite{FR1}.  Choosing $\eta$ small and then $r \in (0,1) $ sufficiently small, we have
\be
rt < \sigma_\Gamma(t) \leq t-1 \no
\ee
whenever $\tau_\eta(t) > t$. Now,
\be
P^{ct,t}(\sigma_h(t) > t, \; \tau_\eta(t) > t) \leq P^{ct,t}( A ) \no
\ee
where $A$ is the set $A = \{ \omega \in \Omega | \;r t < \sigma_\Gamma \leq \sigma_h(t) \}$. This probability is bounded above by
\br
& &P^{ct,t}( A ) \leq \no \\
& & e^{\delta t} E_{ct} \left[ e^{\frac{1}{2} \int_0^{\sigma_\Gamma} \zeta(u(Z,t-s)) ds } e^{-\frac{1}{2} \int_0^{\sigma_\Gamma} \zeta(u(Z,t-s)) ds } (u(Z(\sigma_\Gamma),t - \sigma_\Gamma))^{1/2} \chi_{A}\right]. \no
\er
Here we have used the crucial lower bound (\ref{e24}). The expression on the right is bounded by
\be
e^{\delta t} E_{ct} \left[ e^{\frac{1}{2} \int_0^{\sigma_\Gamma} \zeta(u(Z,t-s)) ds } e^{-\frac{1}{2} r t \zeta_h } (u(Z(\sigma_\Gamma),t - \sigma_\Gamma))^{1/2} \chi_{A}\right]. \no
\ee
Then using H\"older's inequality, we have
\br
P^{ct,t}(\sigma_h(t) > t, \; \tau_\eta(t) > t) &\leq& e^{\delta t} [u(ct,t)]^{1/2} e^{- \frac{1}{2} \zeta_h r t} \no \\
& \leq & e^{t ( \delta - \frac{1}{2} \zeta_h r )}
\er
The last term goes to zero if $\delta$ is sufficiently small, depending on $h$.  This completes the proof of Theorem \ref{theo:frspeed}. \qed

Note that the only difference between this argument and that of Freidlin in the periodic case is the manner in which (\ref{e24}) and (\ref{e16}) are obtained.

\section{Proof of Lemma \ref{lem:explowbnd}}
\setcounter{equation}{0}

The main issue in proving the estimate of Lemma \ref{lem:explowbnd} (and thus the lower bound (\ref{e24})) is whether the random variable
\be
\eta^t_z(\kappa t) = \frac{x - X^{z,t}(\kappa t) }{\kappa t} \label{e106}
\ee
satisfies a large deviation principle with a convex rate function that can be characterized by $\mu(\lambda)$, almost surely with respect to $Q$.  The variable $\eta^t_z(\kappa t)$ is the first component of the average speed of a trajectory over time interval $[0, \kappa t]$.

\begin{defin}\label{def:ldprinciple} For fixed $\hat \omega \in \hat \Omega$, the random variables $\eta^t_z(\kappa t)$ satisfy a large deviation principle with a convex rate function $S(c)$ if there exists a convex function $S(c)$, independent of $z \in R^n$, such that
\bi
\item[(i)] For each $s \geq 0$, the set $\Phi(s)  = \{ c \in R | \; S(c) \leq s \}$ is compact.
\item[(ii)] For any $\delta, h > 0$, there exists $t_0 > 0$ such that for all $t > t_0$
\be
P\left( d(\eta^t_z(\kappa t),\Phi(s)) > \delta \right) \leq e^{-\kappa t(s-h)}. \no
\ee
\item[(iii)] For any $\delta, h > 0$, there exists $t_0 > 0$ such that for all $t > t_0$
\be
P\left( \eta^t_z(\kappa t) \in U_\delta(c) \right) \geq e^{- \kappa t(S(c) + h)}.
\ee
\ei
\end{defin}
If such a function $S(c)$ exists, it might depend on the parameter $\kappa \in (0,1]$, and it might depend on $\hat \omega \in \hat \Omega$. However, we will show that
\begin{prop}\label{prop:rateexist}
Suppose that assumptions A1-A6 hold. Then almost surely with respect to $Q$, the random variables $\eta^t_z(\kappa t)$ satisfy a large deviation principle (with respect to $P$) with a convex rate function $S(c)$ that is independent of $\kappa$ and $\hat \omega \in \hat \Omega$.
\end{prop}
We postpone the proof for the moment while we finish the proof of Lemma \ref{lem:explowbnd}. By Lemma \ref{lem:muexist}, the quantity
\be
\bar \mu(\lambda) = \mu(\lambda) - f'(0) \no
\ee
is well defined and is almost surely constant with respect to $Q$ for $\lambda \in R$, independently of $\kappa$. Since, by our assumption of Proposition \ref{prop:rateexist}, the variables $\eta^t_z(\kappa t)$ have a convex rate function, it follows that $S(c)$ has the following characterization (see Section 5.1 of \cite{FW1}):
\be
S(c) = \sup_{\lambda \in R} [c \lambda -  \bar \mu (\lambda)]. \label{e53}
\ee
This characterization does not hold if $S(c)$ is not convex. Let us emphasize that $S(c)$ is independent of $\kappa \in (0,1]$ and $\hat \omega \in \hat \Omega$, although the constants $t_0$ in Definition \ref{def:ldprinciple} may depend on $\kappa, \hat \omega$.

Now, by definition of $S(c)$,
\be
\liminf_{t \to \infty} \frac{1}{\kappa t} \log \inf_{z \in R^n} P\left\{ \eta^t_z(\kappa t) \in U_\delta(c) \right\} \geq - S(c) > - \infty, \label{e29}
\ee
and the lower bound (\ref{e27}) of Lemma \ref{lem:explowbnd} can be written
\be
\liminf_{t\to \infty} \frac{1}{t} \log \inf_{c \in K, y \in D} u(ct,y,t) \geq  f'(0) - \max_{c \in K} S(c). \label{e28}
\ee
To prove the lower bound we now use the Feynman-Kac formula to relate (\ref{e28}) to (\ref{e29}), as in the arguments of Freidlin in Lemma 7.3.2 in \cite{FR1}. The compactness of $K$ implies that it suffices to show that given any $\epsilon > 0$, and any $c$ for which $H(c) > 0$,
\br
& & \liminf_{t \to \infty} \left ( \frac{1}{t} \log \inf_{\tilde c \in U_\delta(c),\;y \in D} u(\tilde c t,y, t)\right)   \geq f'(0) - S(c) - \epsilon \label{e9}
\er
for $\delta > 0$ sufficiently small. Without loss of generality, we assume that
\be
u_0(x) \geq \chi_{U_\delta(0)}(x) \label{e4}
\ee
for $\delta$ sufficiently small. That is, $u_0 = 1$ whenever, $\abs{x} < \delta$.  We define $q$ to be the limit on the left hand side of (\ref{e9}):
\be
q = \liminf_{t \to \infty} \left ( \frac{1}{t} \log \inf_{x \in U_{\delta t}(ct), y \in D} u(x,y, t)\right) . \no
\ee
Without loss of generality, we assume $0 < c^* < c $ where $c^*$ is defined by (\ref{e65}).

{\bf Step 1:}
The first step is essentially the same as in \cite{FR1}. Suppose for the moment that we known $q$ is finite. By the representation (\ref{fkeq2}) we have for any $\kappa \in (0,1]$
\be
\inf_{\tilde c \in U_\delta(c), y \in D} u(t \tilde c , y, t) \geq  \inf_{\tilde c \in U_\delta(c), y \in D} E \left [e^{\int_0^{\kappa t} \zeta(t-s,u(Z(s),t-s))ds} u(Z(t-\kappa t), t - \kappa t) \chi_A \right] \label{e30}
\ee
for any set $\mathcal{F}_{s \leq t} $-measurable set $A$.  Recall that when $u \leq h$, the reaction rate $\zeta(u)$ is bounded below by $\zeta_h > 0$. If we choose $A$ to be the set of paths satisfying both
\be
X^{z,t}(\kappa t) \in U_{(1-\kappa)\delta t}( (1-\kappa) t c) \label{e104}
\ee
and
\be
u(Z^{z,t}(s),t-s) \leq h \;\;\text{for all}\;\; s \in [0,\kappa t], \label{e31}
\ee
then from (\ref{e30}) and the assumption that $q$ is finite we have a lower bound
\br
q & \geq & \zeta_h  +  \liminf_{t \to \infty} \frac{1}{\kappa t} \log \inf_{\tilde c \in U_\delta(c), y \in D} P (A), \label{e32}
\er
provided that the limit on the right also exists and is finite.

{\bf Step 2:}
Now we bound the right hand side of (\ref{e32}) and show how it relates to (\ref{e29}). Let $\delta$ be sufficiently small so that we can pick $c'$ with  $c^* < c' < c - 6 \delta $. By Theorem \ref{theo:frspeed1}, for any $h \in (0,1)$ there is a constant $t_0> 0$ such that
\be
u(x,y,t) \leq h \;\;\; \text{for all} \;\;x \geq c' t,\; y \in R^{n-1},\; t \geq t_0. \no
\ee
Now if $\kappa < 1/2$ and
\be
\sup_{s \in [0,\kappa t]} \abs{X^{z,t}(s) - (t-s) c } \leq 3 \delta t,
\ee
then $X^{z,t}(s) > c'(t-s) $ for all $s \in [0,\kappa t]$. Thus, (\ref{e31}) is achieved along such paths when $t > 2 t_0$. Next, if $\tilde c \in U_\delta(c)$ is written $\tilde c = c + \delta e_1$ for some $e_1$ with $\abs{e_1} < 1$, then define $\hat c = c + 2 \delta e_1$. Then for any $\abs{e_2} < 1$
\be
\tilde c t - \kappa t \hat c + \kappa t \delta e_2 \in U_{(1-\kappa)\delta t}((1-\kappa)ct).
\ee
It follows that for each $\tilde c \in U_{\delta}(c)$ there is a $\hat c \in U_{2 \delta}(c)$ such that (\ref{e104}) is achieved whenever $\eta^t_z(\kappa t) \in U_\delta(\hat c)$, where $\eta$ is defined by (\ref{e106}). This gives us a lower bound on $P(A)$ in terms of the $\eta^t_z(\kappa t)$, the average speed of a trajectory over $[0,\kappa t]$:
\br
& & \inf_{\tilde c \in U_\delta(c), y \in D, z =  ct} P (A)  \geq  \label{e105} \\
 & & \inf_{\hat c \in U_{2 \delta}(c), y \in D, z = \hat ct} P \left(\sup_{s \in [0,\kappa t]} \abs{X^{z,t}(s) - (t-s) c } \leq 3 \delta t , \;\; \eta^t_z(\kappa t) \in U_\delta(\hat c) \right) \no
\er

For $\kappa $ sufficiently small, $\kappa <  (2 \delta)/(3 \max(1,\abs{c}))$, we see that
\be
\sup_{\hat c \in U_{2 \delta}(c), y \in D, z = \hat ct} P\left (\sup_{s \in [0,\kappa t]} \abs{X^{z,t}(s) - (t-s) c } \geq 3 \delta t  \right) \leq  \sup_{z \in R^n} P \left (\sup_{s \in [0,\kappa t]} \abs{X^{z, t}(s) - x} \geq  \delta t/3 \right ). \no
\ee
By Lemma \ref{corol:taubound} there are constants $K_1, K_2 > 0$ independent of $\kappa$ such that (except possibly on a set of $Q$-measure zero)
\be
\sup_{z \in R^n} P \left (\sup_{s \in [0,\kappa t]} \abs{X^{z, t}(s) - x} \geq  \delta t/3 \right ) \leq K_1 e^{-K_2 \delta^2 t /\kappa}  \label{e45}
\ee
for $t$ sufficiently large depending on $\hat \omega$.  Therefore, for any $M > 0$, by choosing $\kappa$ arbitrarily small, we can make $K_2 \delta^2/\kappa^2 > M$, so that
\br
& & \limsup_{t \to \infty} \frac{1}{\kappa t} \log(\sup_{z \in R^n} P\left( \sup_{s \in [0,\kappa t]} \abs{X^{z,t}(s) - (t-s) c } \geq 2 \delta t  \right)) \leq  \no \\
& & \leq \lim_{t \to \infty} \frac{1}{\kappa t} \log( K_1 e^{-K_2 \delta^2 t /\kappa} )\leq   -M . \no
\er
Therefore, from (\ref{e105}) and (\ref{e32}) we now see that for $\kappa$ sufficiently small,
\br
q & \geq & \zeta_h  +  \liminf_{t \to \infty} \frac{1}{\kappa t} \inf_{\hat c \in U_{2 \delta}(c), z \in R^n} P \left(\eta^t_z(\kappa t) \in U_\delta(\hat c) \right)
\er
provided that the limit on the right is finite and bounded below, independently of $\kappa$. However, this follows immediately from Proposition \ref{prop:rateexist} and the lower bound (\ref{e29}). Then (\ref{e9}) follows by letting $h \to 0$ so that $\zeta_h \to f'(0)$.

{\bf Step 3:}
It remains to establish the initial claim that $q > -\infty$. To see this, define for any $c_1 \in R$
\be
\hat q_\delta(c_1,t) = \inf_{x \in U_{\delta}(tc_1),\, y \in D} P_x \left( X^{z,t}(t) \in U_{\delta}(0) \right), \label{e5}
\ee
a random variable over $\hat \Omega$.  Let use write $X^{x,y,t}(s)$ as
\be
X^{z,t}(s) = x + I^{y,t}(s) + W_1(s) \no
\ee
where $I^{y,t}$ is the first integral term in (\ref{e12}) and $W_1(0) = 0$. For $\hat \omega \in \Omega$ fixed, let $M > 0$ and define the set
\be
A_M = \{ w |\; \sup_{y \in D, s \in [0,t]} \abs{I^y(s)} \leq M t \} . \no
\ee
Using the fact that $W_1$ and $W_2$ are independent, we derive the lower bound for $x \in U_\delta(c_1 t)$
\br
P(X^{z,t}_t \in U_\delta(0)) & \geq &  P(W_1(t) \in U_\delta(0) - I^y(t) - x ;\; A) \no \\
& \geq & P(W_1(t) \in U_\delta(0) + Mt + \abs{c_1} t + \delta)P(A_M) \no \\
& \geq & \frac{\delta}{\sqrt{2 \pi t}} e^{- \frac{(Mt + \abs{c_1}t + 2\delta)^2}{ 2 t}}P(A_M).  \label{e68}
\er
By Lemma \ref{corol:taubound}, $P(A_M) \geq 1/2$ for $t$ sufficiently large, depending on $\hat \omega$ and $M$. Therefore, there is a finite constant $K_1 > 0$ depending only on $c_1$ and $M$ such that
\be
\liminf_{t \to \infty} \frac{1}{t} \log \hat q_\delta(c_1,t) = -K_1 \label{e71}
\ee
is finite almost surely with respect to $Q$. Now for $c_1 < c_2$, we get two constants $K_1$ and $K_2$ such that the limit (\ref{e71}) holds almost surely with respect to $Q$. Let $c_1 = c - \delta$ and $c_2 = c+\delta$ and define the Markov time
\be
\pi(z,t) = \inf_s \left\{ s > 0 | \; X^{z,t}(s) \in (-\infty, c_1 (t-s)] \cup [c_2 (t-s) , \infty);\;\; \text{or}\; s \geq t \right \} \no
\ee
which is the first exit time of the process $ X^{x,t}(s) $ from the region defined by $\{ (x,y,\tau) | \; x \in (c_1 \tau , c_2 \tau) , \tau \geq 0, y \in D\}$.
By formula (\ref{fkeq}), we have
\br
\inf_{\tilde c \in U_\delta(c), y \in D} u(\tilde c t, y, t) & = & \inf_{x \in (c_1 t, c_2 t), y \in D} u(x, y, t)  \label{e8} \\
& \geq & \inf_{x \in (c_1 t, c_2 t), y \in D} E \left [ u(Z^{z,t}(\pi), t - \pi)\right] \no \\
& \geq & \left( \inf_{s  \in [0,t] , y \in D} u(c_1 s,y,s)\right) \wedge  \left(\inf_{s  \in [0,t] , y \in D} u(c_2 s,y,s)\right). \no
\er
It follows from (\ref{e4}) and (\ref{e71}), that there are constants $C_1, C_2 > 0$ (depending on $\hat \omega$) such that
\br
u(c_1 s,y,s)  \geq   C_1 e^{- 2 K_1 s}  \;\;\; \text{and} \;\;\; u(c_2 s,y,s)  \geq   C_2 e^{- 2 K_2 s} \label{e7}
\er
for all $s \geq 0$. Combining (\ref{e8}) with (\ref{e7}), it is clear that $q \geq - 2 (K_1 \vee K_2) > - \infty$. Having shown that $q$ is finite, we have completed the proof of Lemma \ref{lem:explowbnd}. \qed

For later use, let us show that for all $t>0$, $\log(\hat q_\delta(c,t))$ is integrable with respect to $Q$. Using (\ref{e68}) we see that
\br
\frac{1}{t} \log \hat q_\delta(c,t) & \geq &  \frac{1}{t} \log \left( \frac{\delta}{\sqrt{2 \pi t}} e^{- \frac{(Mt + \abs{c_1}t + 2\delta)^2}{ 2 t}}P(A_M)\right) \\ \no
& \geq &  - C_1 + \frac{1}{t} \log \left( e^{- \frac{(Mt + \abs{c_1}t + 2\delta)^2}{ 2 t}}P(A_M)\right)  \no
\er
for a constant $C_1>0$ independent of $c$, for $t \geq 1$. Let $\hat g$ be the term inside the logarithm:
\be
\hat g = e^{- \frac{(Mt + \abs{c_1}t + 2\delta)^2}{ 2 t}}P(A_M).
\ee
Then for $\alpha \geq 2 C_1$,
\br
Q\left( \frac{1}{t} \log \hat q_\delta(c,t) \leq - \alpha\right) &\leq& Q\left( \frac{1}{t} \log \hat g \leq - \alpha/2 \right) \no \\
&=& Q\left(\hat g \leq e^{- \alpha t/2} \right) \no \\
&=& Q\left( P(A_M) \leq e^{- \alpha t/2} e^{\frac{(Mt + \abs{c_1}t + 2\delta)^2}{ 2 t}} \right). \label{e76}
\er
Also, from Lemma \ref{lem:borelbound},
\be
Q \left( P(A_M) \leq 1 - e^{-K_2 M^2 t/2} \right) \leq K_1 e^{-K_2 M^2 t/2}. \label{e77}
\ee
With a little algebra, one can see that there exist constant $K_3, K_4 > 0$ independent of $t$ such that whenever $t > 1$, $M = K_3 \sqrt{\alpha}$, and $\alpha \geq K_4 c^2$ , we have
\be
 e^{- \alpha t/2} e^{\frac{(Mt + \abs{c_1}t + 2\delta)^2}{ 2 t}} \leq 1/2 \leq 1 - e^{-K_2 M^2 t/2} . \no
\ee
By combining (\ref{e76}) and (\ref{e77}), we now conclude that
\br
Q\left( \frac{1}{t} \log \hat q_\delta(c,t) \leq - \alpha\right) &\leq& K_1 e^{-K_2 K_3^2 \alpha t/2} \label{e79}
\er
whenever $\alpha \geq K_4 c^2$ and $t > 1$. It follows that for $t$ sufficiently large, independent of $\hat \omega$,
\br
E_Q[ \abs{\frac{1}{t} \log \hat q_\delta(c,t)} ]  &=& \int_0^\infty Q\left( \abs{\frac{1}{t} \log \hat q_\delta(c,t)} \geq \alpha \right)\, d\alpha  \no \\
& \leq & K_4 c^2 + \int_{K_4 c^2}^\infty  K_1 e^{-K_2 K_3^2 \alpha t/2} \, d\alpha < \infty . \label{e80}
\er
This is bounded uniformly in $t$, for $t > 1$.

\section{Proof of Large Deviation Estimates}
In this section we prove Proposition \ref{prop:rateexist}. We work first with the case $\kappa = 1$. Because $b$ is independent of $x$, it suffices to show that the proposition holds for $x = 0$. Define for $c \in R$ and $r < s < t$
\be
q^y_\delta(c,r,t) =  P(-X^{0,y,t}(t-r) \in U_{\delta (t-r)}(c(t-r))) = P(\eta^t_y(t-r) \in U_\delta(c))\no
\ee
and
\br
q^+_\delta(c,r,t) & = & \sup_y q^y_\delta(c,r,t)  \no \\
q^-_\delta(c,r,t) & = & \inf_y q^y_\delta(c,r,t).  \no
\er
We will use the subadditive ergodic theorem to show that $(1/t) \log q^-_\delta(c,0,t)$ and $(1/t) \log q^+_\delta(c,0,t)$ converge almost surely to a finite constant. Define the events
\br
A & = & \left\{ -X^{0,y,t}(t-r) \in U_{\delta (t-r)}(c(t-r)) \right\} = \left\{ \eta^t_y(t-r) \in U_\delta(c) \right\} \no \\
B & = & \left\{ -X^{0,y,t}(t-s) \in U_{\delta (t-s)}(c(t-s)) \right\} = \left\{ \eta^t_y(t-s) \in U_\delta(c) \right\}. \no
\er
Note that event $B$ is $\mathcal{F}_{s \leq \tau}$ measurable for $\tau \geq t-s$, where $\mathcal{F}_\tau$ is the $\sigma$-algebra generated by $(Z^{z,t}(s))_{s \leq \tau}$.  By the Markov property of the Wiener process we have:
\br
q^-_\delta(c,0,t) & = & \inf_y P(A) \geq \inf_y P(A \cap B) \no \\
& \geq & \inf_y P( \left\{-X^{0,y,t}(t-r) + X^{0,y,t}(t-s) \in U_{\delta (r-s)}(c(s-r))\right\} \cap B) \no \\
& = & \inf_y E\left[ \chi_B P[ \left\{-X^{0,y,t}(t-r) + X^{0,y,t}(t-s) \in U_{\delta (s-r)}(c(s-r))\right\}  \;|\;\mathcal{F}_{t-s} ] \right] \no \\
& = & \inf_y E\left[ \chi_B P[ \left\{-X^{0,y,t}(t-r) + X^{0,y,t}(t-s) \in U_{\delta (s-r)}(c(s-r))\right\}  \;|\; X^{0,y,t}(t-s) ] \right] \no \\
& \geq & \inf_y E\left[ \chi_B \inf_y P[ \left\{-X^{0,z,s}(s-r) \in U_{\delta (s-r)}(c(s-r)) \right\}  ] \right] \no \\
& = &  \inf_y P( \left\{-X^{0,z,s}(s-r) \in U_{\delta (s-r)}(c(s-r)) \right\} ) \inf_y P ( B ) \no \\
& = & q^-_\delta(c,r,s)q^-_\delta(c,s,t). \no
\er
Therefore, $\log(q^-_\delta(c,s,t))$ is super-additive for each $c \in R$. By the stationarity of $b$,
\br
\tau_h q^-_\delta(c,r,t) &  = & \tau_h \inf_y P(-X^{0,y,t}(t-r) \in U_{\delta (t-r)}(c(t-r))) \no \\
& = & \inf_y P(-X^{0,y,t+h}(t - r) \in U_{\delta (t-r)}(c(t-r))) \no \\
& = & q^-_\delta(c,r+h,t+h). \no
\er
For any $\epsilon > 0$, we can bound $q$ below by translating in $x$ and using (\ref{e5}):
\be
q^-_\delta(c,r,t) \geq \tau_{r}\hat q_\epsilon(c,t-r) = \inf_{y \in D, x \in U_\epsilon(ct)} P(X^{x,y,t}(t-r) \in U_{\epsilon}(cr)) \label{e69}
\ee
for $t-r$ sufficiently large. Hence, $\log(q^-_\delta(c,r,t))$ is integrable by (\ref{e80}). Kingman's ergodic theorem \cite{King} now implies that the limit
\be
\lim_{n \to \infty} \frac{1}{n} \log q^-_\delta(c,0,n) = \sup_{n > 0} \frac{1}{n} E_Q [ \log q^-_\delta(c,0,n) ] = - S_\delta(c) \label{e52}
\ee
exists and is a finite constant, $Q$-a.s, because of the ergodicity assumption A3. Using an idea in \cite{AK}, we can now extend convergence in (\ref{e52}) to continuous time. Let
\be
g(\hat \omega) = \sup_{\substack{r,t \in [0,2]\\ \abs{r-t} \geq 1 }} \abs{\log(q^-_\delta(c,r,t))} . \no
\ee
Let $\Upsilon(\hat \omega) = \sup_{y \in D, t \in [0,2]} \abs{b(y,t)}$.
Then for all $0 \leq r < t \leq 2$,
\be
 \sup_{y \in D} \abs{ \int_0^{t-r} b(W^y(\tau), t- \tau ) d\tau \,}\leq \Upsilon \abs{t-r}. \no
\ee
As in (\ref{e68}),
\br
P(-X^{0,y,t}(t-r) \in U_{\delta (t-r)}(c(t-r))) & \geq & P(W_1(t-r) \in U_{\delta(t-r)}(0) + \Upsilon (t-r) + \abs{c}(t-r)) \no \\
& \geq & \frac{\delta \abs{t-r}}{\sqrt{2 \pi \abs{t-r}}} e^{- \frac{(t-r)^2(\Upsilon  + \abs{c} + \delta)^2}{ 2 (t-r)}}.   \no
\er
Therefore,
\be
0 \geq g(\hat \omega) \geq - K_1  - K_2 \Upsilon^2 \no
\ee
for some constants $K_1,K_2>0$ that depend on $\delta$ and $c$. Hence $g(\hat \omega)$ is integrable with respect to $Q$. By the super-additivity of $q^-_\delta(c,r,t)$,
\be
q^-_\delta(c,0,n-1) - \tau_{n-1} g \leq q^-_\delta(c,0,t) \leq q^-_\delta(c,0,n+2) + \tau_{n} g \label{e113}
\ee
whenever $t \in (n ,n+1)$, $n \in Z$.  As in the proof of Theorem 2.5 of \cite{AK}, one can show that $\frac{1}{n} \tau_n g \to 0$ almost surely as $n \to \infty$ since $g$ is integrable. It now follows from (\ref{e113}) that
\be
\lim_{t \to \infty} \frac{1}{t} \log q^-_\delta(c,0,t) = - S_\delta(c) \no
\ee
almost surely with respect to $Q$.

For each $c \in R$, $S_\delta(c)$ can be bounded above independently of $\delta > 0$ using (\ref{e69}) and (\ref{e71}). From the definition, it is clear that $S_\delta(c) \geq 0$ for all $\delta$, and $S_{\delta_1}(c) \leq S_{\delta_2}(c)$ whenever $\delta_1 > \delta_2$.  Therefore, we define for each $c \in R$
\be
S(c) = \lim_{\delta \to 0} S_\delta(c) = \sup_{\delta > 0} S_\delta(c) \in [0,+\infty). \no
\ee
\begin{lem}
For all $\delta > 0$, the functions $S_\delta(c)$ are continuous and convex in $c$. Also, $S(c)$ is continuous and convex in $c$.
\end{lem}
{\bf Proof:}
The continuity and convexity of $S(c)$ follows immediately from the continuity and convexity of $S_\delta(c)$, since the functions $S_\delta(c)$ converge pointwise to the finite limit $S(c)$. The convexity of $S_\delta(c)$ follows from the Markov property of the process $X^{z,t}$, as follows.

Let $p \in (0,1)$ and $c_0 = p c_1 + (1-p) c_2$. Let $t > 0$ and denote $t_1 = p t$, $t_2 = (1-p) t$. Then we see that
\br
q^-_\delta(c_0,0,t) & = & \inf_y P\left(-X^{0,y,t}(t) \in U_{\delta t}(c_0 t)\right)\no \\
& \geq & \inf_y P\left(-X^{0,y,t}(t) \in U_{\delta t}(c_0 t) \;,\;\; -X^{0,y,t}(t_1) \in U_{\delta t_1}(c_1 t_1) \right) \no \\
& \geq & \inf_y P\left(-X^{0,y,t-t_1}(t_2) \in U_{\delta t_2}(c_2 t_2)\right ) \inf_y P\left(-X^{0,y,t}(t_1) \in U_{\delta  t_1 }(c_1 t_1) \right) \no \\
& = & q^-_\delta(c_2,0,t_2)  q^-_\delta(c_1,t_2, t) . \no
\er
Hence
\br
-\frac{1}{t} \log q^-_\delta(c_0,0,t) & \leq &  - \frac{1}{t} \log q^-_\delta(c_2,0,(1-p) t) - \frac{1}{t} \log q^-_\delta(c_1,(1-p)t,t) \label{e57} \\
& = &  - \frac{1}{t} \log q^-_\delta(c_2,0,(1-p) t) -  \tau_{(1-p)t} \left(\frac{1}{t} \log q^-_\delta(c_1,0,pt)\right).\no
\er
By the stationarity of $b$, the random variable
\be
-\tau_{(1-p)t} \left(\frac{1}{t} \log q^-_\delta(c_1,0,pt)\right) \no
\ee
converges in distribution to $p S_\delta(c_1)$. Therefore, there is a set $M \subset \hat \Omega$ with $Q(M) = 0$ such that for each $\hat \omega \in \hat \Omega \setminus M$ and each $\epsilon > 0$, there is an increasing sequence $\{t_j \}_{j=1}^\infty$, $t_j \to \infty$ as $j \to \infty$, such that for $j $ sufficiently large,
\be
- \tau_{(1-p)t_j} \left(\frac{1}{t_j} \log q^-_\delta(c_1,0,pt_j)\right) \leq p S_\delta(c_1) + \epsilon. \label{e54}
\ee
Note that the other two terms in (\ref{e57}) are random variables that converge to constants ($Q$-a.s.) as $t \to \infty$. For fixed $\hat \omega \in \hat \Omega$ and $\epsilon < 0$, we can pick $t$ sufficiently large so that
\be
- \frac{1}{t} \log q^-_\delta(c_2,0,(1-p) t) \leq (1-p) S_{\delta}(c_2) + \epsilon, \label{e55}
\ee
and
\be
-\frac{1}{t} \log q^-_\delta(c_0,0,t) \geq S_\delta(c_0) - \epsilon \label{e56}.
\ee
Now by (\ref{e57}) and (\ref{e54})-(\ref{e56}) we have
\be
S_\delta(c_0) - \epsilon \leq (1-p) S_{\delta}(c_2) + p S_\delta(c_1)  + 2 \epsilon. \no
\ee
Since $\epsilon >0$ was arbitrarily chosen, we infer that
\be
S_\delta(c_0) \leq (1-p) S_{\delta}(c_2) + p S_\delta(c_1) .
\ee
So, $S_\delta(c)$ is convex and must also be continuous in $c$, since it is finite for every $c \in R$. \qed

This establishes the existence and convexity of the function $S(c)$. Part (iii) of the Definition \ref{def:ldprinciple} follows from the definition of $S_{\delta}(c)$ and the fact that $S_{\delta}(c) \nearrow S(c)$.

To finish the proof of Proposition \ref{prop:rateexist} for $\kappa =1$, we must establish a Harnack-type inequality to relate the quantities
\be
P(-X^{0,y,t}(t) \in U_{\delta t}(ct)) \;\;\;\text{and} \;\;\; P(-X^{0,y',t}(t) \in U_{\delta t}(ct)) \no
\ee
for $y, y' \in D$. This will allow us to remove the $\inf_y$ in the definition of $q$ and $S(c)$ and to establish parts (i) and (ii) of Definition \ref{def:ldprinciple}. We prove the following lemma
\begin{lem} \label{lem:harnackP}
There are constants $K_1, K_2, K_3 > 0$ such that for all $\kappa \in (0,1]$, $c \in R$, $\epsilon>0$, and $\delta > 0$,
\be
\inf_z P \left(\eta^t_z(\kappa t) \in U_{(1 + \epsilon) \delta }(c) \right) \geq K_3 \sup_z P \left(\eta^t_z(\kappa t) \in U_{ \delta }(c) \right) - K_1 e^{-K_3 \epsilon^2 \delta^2 \kappa^2 t^2} \no
\ee
\end{lem}
{\bf Proof:}
For clarity we let $\kappa = 1$. Extension to $\kappa < 1$ is straightforward, as in the proof of Lemma \ref{corol:taubound} in the appendix. First, note that
\be
X^{x,y,t}(t) = x_0 + I^{y,t}(t) + W_1(t), \;\;W^1_0 = 0 \no
\ee
and
\be
X^{x,y,t}(t) - X^{x,y,t}(s) = I^{y + W_2(s) , s}(t-s) + W_1(t) - W_1(s). \no
\ee
With out loss of generality, we assume $x_0 = 0$.  Therefore,
\br
& & P(X^{0,y,t}(t) \in U_{\delta t}(-ct)) \no \\
& & \;\;\;\;\;=  P(X^{0,y,t}(t) - X^{0,y,t}(s)   \in U_{\delta t}(-ct - X^{0,y,t}(s) )) \no \\
& & \;\;\;\;\;=  E\left( P( I^{\beta , s}(t-s) + W_1(t) - W_1(s) \in U_{\delta t}(-ct - \alpha ) | X^{0,y,t}(s) = \alpha , \, W_2(s) = \beta - y)\right) \no \\
& & \;\;\;\;\;=  E\left( P( I^{\beta , s}(t-s) +  \tilde W_1(t-s) \in U_{\delta t}(-ct - \alpha ) | X^{0,y,t}(s) = \alpha , \, W_2(s) = \beta - y)\right). \no
\er
Here we used $\tilde W_1(\tau)$ to denote a Wiener process, $\tilde W_1(0) = 0$, P-a.s. Now we see that
\br
& & P(X^{0,y,t}(t) \in U_{\delta t}(-ct), \;\; \abs{ X^{0,y,t}(s)} \leq M) \no \\
& & \;\;\; \leq  E\left( P( I^{\beta , t-s}(t-s) + \tilde W_1(t-s) \in U_{\delta t + M}(-ct ) | \, W_2(s) = \beta - y)\right) \no
\er
and
\br
& & P(X^{0,y,t}(s) \in U_{\delta t}(-ct), \;\; \abs{ X^{0,y,t}(s)} \leq M) \no \\
& & \;\;\; \geq  E\left( P( I^{\beta , t-s}(t-s) + \tilde W_1(t-s) \in U_{\delta t - M}(-ct ) | \, W_2(s) = \beta - y)\right). \no
\er
For any $y, y' \in D$, we have for fixed $s \in (0,t]$
\br
& & P(X^{0,y',t}(t) \in U_{\delta t + 2M}(-ct), \;\; \abs{ X^{0,y',t}(s)} \leq M) \no \\
& & \;\;\; \geq  E\left( P( I^{\beta , t-s}(t-s) + \tilde W_1(t-s) \in U_{\delta t + M}(-ct ) | \, W_2(s) = \beta - y')\right) \no \\
& & \;\;\; = \int_D \rho_s(\beta - y') P( I^{\beta , t-s}(t-s) + \tilde W_1(t-s) \in U_{\delta t + M}(-ct )) \,d\beta \no \\
& & \;\;\; \geq C_1 \int_D \rho_s(\beta - y) P( I^{\beta , t-s}(t-s) + \tilde W_1(t-s) \in U_{\delta t + M}(-ct )) \,d\beta \no \\
& & \;\;\; = C_1 E\left( P( I^{\beta , t-s}(t-s) + \tilde W_1(t-s) \in U_{\delta t + M}(-ct ) | \, W_2(s) = \beta - y)\right) \no \\
& & \;\;\; \geq C_1 P(X^{0,y,t}(t) \in U_{\delta t }(-ct), \;\; \abs{ X^{0,y,t}(s)} \leq M)
\er
where $\rho_s(r)$ denotes density function for distribution of $W_2(s)$ on the torus, and $C_1$ is a positive constant depending on $s$. Using this inequality, we see that
\br
P(X^{0,y',t}(t) \in U_{\delta t + 2M}(-ct))  & = & P(X^{0,y',t}(t) \in U_{\delta t + 2M}(-ct), \;\; \abs{ X^{0,y',t}(s)} \leq M)  \no \\
 & & + \;P(X^{0,y',t}(t) \in U_{\delta t + 2M}(-ct), \;\; \abs{ X^{0,y',t}(s)} > M) \no \\
& \geq & C_1 P(X^{0,y,t}(t) \in U_{\delta t }(-ct))  \no \\
&  & - \; C_1 P(X^{0,y,t}(t) \in U_{\delta t }(-ct), \;\; \abs{ X^{0,y,t}(s)} > M)  \no \\
 & & + \;P(X^{0,y',t}(t) \in U_{\delta t + 2M}(-ct), \;\; \abs{ X^{0,y',t}(s)} > M). \no \\
 \label{e47}
\er
For $\epsilon > 0$, let $s=1$ and $M = \epsilon \delta t/2$ ($M = \epsilon \delta \kappa t/2$ when $\kappa < 1$). It follows from the Borell inequality, as in proof of Lemma \ref{corol:taubound}, that
\be
P\left( \sup_{\tau \in [0,1], y \in D} \abs{X^{0,y,t}(\tau)} > M  \right) \leq K_1 e^{-K_2 M^2} =  K_1 e^{-K_2 \epsilon^2 \delta^2 t^2 / 4}
\ee
for $t$ sufficiently large. Now the lemma follows from (\ref{e47}). \qed

Since $\epsilon > 0$ is arbitrary, the lemma implies that
\begin{cor}
\br
S_\delta(c) &=& \lim_{t \to \infty} \frac{1}{t} \log q^-_\delta(c,0,t) =  \lim_{t \to \infty} \frac{1}{t} \log q^+_\delta(c,0,t) \label{e50}
\er
almost surely with respect to $Q$.
\end{cor}
Now, using this estimate, we can establish parts (i) and (ii) from Definition \ref{def:ldprinciple}. From Lemma \ref{corol:taubound} there are constants $K_1, K_2 > 0$ such that for $t$ sufficiently large,
\be
P( \abs{\eta^t_z(t)} \geq c) \leq K_1 e^{-K_2 c^2 t}. \label{e48}
\ee
This implies that $\lim_{\abs{c} \to \infty} S(c)/\abs{c} = + \infty$. Hence, $\Phi(s)$ is a bounded set, for each $s \geq 0$. Then by continuity of $S(c)$, $\Phi(s)$ is compact. Let $A$ be the set
\be
A = \{ c \in R | \; d(\eta^t_z(t),\Phi(s)) > \delta \}. \no
\ee
We must show that for any fixed $\delta > 0$, $h > 0$, there is a $t_0 > 0$ such that for $t \geq t_0$,
\be
P(\eta^t_z(t) \in A) \leq e^{-t(s - h)}. \label{e49}
\ee
Because of the bound (\ref{e48}), it suffices to show that (\ref{e49}) holds with $A$ replaced by any compact subset $A'$ of $A$ (because $K_2 c^2 > s$ when $c$ is sufficiently large). Pick $\epsilon >0$ small enough (at most $\epsilon < \delta $) so that
\be
\inf_{c' \in A'} S_\epsilon(c') > s - h . \label{e51}
\ee
This is possible since $A'$ is compact and $S_\delta(c') \nearrow S(c')$ as $\delta \to 0$, for any $c' \in A'$. Cover $A'$ with a finite number of $\epsilon$-balls. Because the number of balls is finite, we conclude from (\ref{e50}) and (\ref{e51}) that for $t$ sufficiently large
\be
P(\eta \in A') \leq e^{- t (s - h) } . \no
\ee
This establishes Proposition \ref{prop:rateexist} in the case that $\kappa = 1$. Now we extend the result to $\kappa < 1$, as well. If $\kappa \in (0,1)$ and $\delta > 0$, the stationarity of $b(y,t)$ implies that for any $z \in R^n$
\be
-\frac{1}{\kappa t} \log P\left(\eta^t_z(\kappa t) \in U_\delta(c) \right) = -\frac{1}{\kappa t} \log q^y_\delta((1-\kappa)t,t) \to  S_\delta(c) \label{e58}
\ee
in distribution (with respect to $Q$) as $t \to \infty$, but this does not imply pointwise convergence.

By definition, $q^-_\delta(s,t) \leq q^y_\delta(s,t) \leq q^+_\delta(s,t)$ for all $y,s,t$. Also, as we have already shown, the Markov property implies sub and superadditivity:
\be
\log q^-_\delta(0,t) \geq \log q^-_\delta(0,r) + \log q^-_\delta(r,t) \no
\ee
and
\be
\log q^+_\delta(0,t) \leq \log q^+_\delta(0,r) + \log q^+_\delta(r,t) \no
\ee
for all $r,t > 0$. In the same way, we also have
\be
\log q^y_\delta(0,t) \geq \log q^-_\delta(0,r) + \log q^y_\delta(r,t) \label{eq1}
\ee
and
\be
\log q^y_\delta(0,t) \leq \log q^+_\delta(0,r) + \log q^y_\delta(r,t) \label{eq2}
\ee
for all $y \in D$, and all $r,t > 0$.
For fixed $\kappa$, let $r = (1-\kappa)t$. Then plugging into (\ref{eq1}) we have
\be
\frac{1}{t}\log q^y_\delta(0,t) \geq \frac{1-\kappa}{(1-\kappa) t} \log q^-_\delta(0,(1-\kappa)t) + \frac{\kappa}{\kappa t}\log q^y_\delta((1-\kappa)t,t) . \label{eq6}
\ee
We have already shown that
\be
\lim_{t \to \infty} \frac{1}{t} \log q^-_\delta(0,t) = \lim_{t \to \infty} \frac{1}{t} \log q^+_\delta(0,t) = -S_\delta(c). \no
\ee
It follows from (\ref{eq6}) that
\br
-S_\delta(c) & = &\limsup_{t \to \infty} \frac{1}{t}\log q^y_\delta(0,t) \no \\
& \geq & \lim_{t \to \infty} \frac{1-\kappa}{(1-\kappa) t} \log q^-_\delta(0,(1-\kappa)t) + \limsup_{t \to \infty} \frac{\kappa}{\kappa t}\log q^y_\delta((1-\kappa)t,t) \no \\
& = & - (1-\kappa ) S_\delta(c) + \kappa \limsup_{t \to \infty} \frac{1}{\kappa t}\log q^y_\delta((1-\kappa)t,t). \no
\er
Thus
\be
\limsup_{t \to \infty} \frac{1}{\kappa t}\log q^y_\delta((1-\kappa)t,t) \leq - S_\delta(c). \label{eq3}
\ee
An analogous argument with (\ref{eq2}) shows that
\be
\liminf_{t \to \infty} \frac{1}{\kappa t}\log q^y_\delta((1-\kappa)t,t) \geq - S_\delta(c). \label{eq4}
\ee
Therefore, for each $y \in D$,
\be
\lim_{t \to \infty} \frac{1}{\kappa t}\log q^y_\delta((1-\kappa)t,t) = -S_\delta(c) \label{eq7}
\ee
almost surely with respect to $Q$. Now we need to show that for each $\hat \omega$, the limit is uniform in $y \in D$.
By Lemma \ref{lem:harnackP}, there are constants $K_1$, $K_2$, and $K_3$ such that
\be
q^+_\delta((1-\kappa)t, t) \leq K_1 q^-_{(1+\epsilon)\delta}((1-\kappa)t, t) + K_2 e^{-K_3 \epsilon^2 \delta^2 \kappa^2 t^2}. \no
\ee
It follows from (\ref{eq7}) that for any $y_0$ fixed and any $\epsilon > 0$ and $\delta' = \delta/(1 + \epsilon)$, we have
\br
\liminf_{t \to \infty} \frac{1}{\kappa t}\log q^-_{\delta}((1-\kappa)t,t)  & \geq & \liminf_{t \to \infty} \frac{1}{\kappa t}\log \left( K_1 q^+_{\delta'}((1-\kappa)t, t) - K_2 e^{-K_3 \epsilon^2 \delta^2 \kappa^2 t^2} \right) \no \\
& \geq &  \liminf_{t \to \infty} \frac{1}{\kappa t}\log \left( K_1 q^{y_0}_{\delta'}((1-\kappa)t, t) - K_2 e^{-K_3 \epsilon^2 \delta^2 \kappa^2 t^2} \right) \no \\
& \geq &  - S_{ \delta'}(c) = -S_{\delta/(1 + \epsilon)}(c). \label{eq8}
\er
Since $S_\delta(c) \nearrow S(c)$, inequality (\ref{eq8}) implies that for any $\delta > 0$ and $h > 0$,
\be
\inf_{z \in R^n} P(\eta^t_z(\kappa t) \in U_\delta(c))  \geq e^{-\kappa t (S(c)  + h )} \no
\ee
for $t$ sufficiently large. This proves part (iii) of Definition \ref{def:ldprinciple} for $\kappa \in (0,1]$. Using Lemma \ref{lem:harnackP} and (\ref{eq7}) we also have
\br
\limsup_{t \to \infty} \frac{1}{\kappa t}\log q^+_\delta((1-\kappa)t,t)  & \leq & \limsup_{t \to \infty} \frac{1}{\kappa t}\log \left( K_1 q^-_{(1+\epsilon)\delta}((1-\kappa)t, t) + K_2 e^{-K_3 \epsilon^2 \delta^2 \kappa^2 t^2} \right) \no \\
& \leq &  \limsup_{t \to \infty} \frac{1}{\kappa t}\log \left( K_1 q^{y_0}_{(1+\epsilon)\delta}((1-\kappa)t, t) + K_2 e^{-K_3 \epsilon^2 \delta^2 \kappa^2 t^2} \right) \no \\
& \leq &  - S_{(1 + \epsilon) \delta}(c) . \label{eq5}
\er
Then letting $\delta \to 0$,  $S_{(1+\epsilon) \delta}(c) \nearrow S(c)$, so that (\ref{eq7}) and (\ref{eq5}) imply
\be
P(\eta^t_z(\kappa t) \in U_\delta(c))  \leq e^{-\kappa t (S(c) - h )} \no
\ee
for $t$ sufficiently large and $\delta$ sufficiently small. Now using Lemma \ref{corol:taubound} as in the case of $\kappa = 1$, part (ii) of Definition \ref{def:ldprinciple} follows. This completes the proof of Proposition \ref{prop:rateexist}. \qed

\section{Computing with Variational Formula}
\setcounter{equation}{0}
In this section, we use formula (\ref{e65}) to compute the propagation speed $c^*$.
We also derive some analytical bounds on $c^*$ and compare them with our numerical results.
In our numerical simulations we consider a specific case of $b(y,t) = b_1(y)b_2(t)$, $y \in R^1$, where $b_1(y)$ is a smooth periodic function of $y$ and $b_2(t)$ is stationary Gaussian field. For $b_2$ we use the Ornstein-Uhlenbeck process which solves the It\^o equation
\be
db(t) = -a\, b(t)\,  dt + \, r \, dW(t) , \;\;\; t \geq 0, \label{e43}
\ee
where $W(y,\omega)$ is the standard Wiener process, $X(0,\omega) = X_0(\omega)$ is a Gaussian random variable with mean zero, and variance $\rho = r^2/(2a)$. This is a mean zero stationary Gaussian process with covariance function
\be
E_Q[b_2(t)b_2(s)] = \Gamma(\abs{t-s}) = \frac{r^2}{2a}e^{- a\abs{t-s}}. \no
\ee
In most computations, we will choose $r = \sqrt{2} \alpha^{3/4}$ so that the covariance is
\be
E_Q[b_2(t)b_2(s)] = \sqrt{\alpha} e^{-\alpha \, |t-s|}  = V(\abs{t-s}) \quad. \label{e91}
\ee
By this choice of $r$, the $L^2$ norm of $V(z)$ remains constant as $\alpha$ changes, so that the total energy in the power spectrum of the $b_2$ remains constant. Because $b_1(y)$ is periodic and because of the rapid decay of the covariance function $V$, it is easy to verify that $b(y,t) = b_1(y)b_2(t)$ satisfies the assumptions A1-A5. Assumption A6 follows from the Lipschitz continuity of $b_1(y)$.

\subsection{Numerical Computation of $\mu(\lambda)$}
To compute $\mu(\lambda) = \lambda^2/2 + f'(0) + \rho(\lambda)$ (see \ref{asle2}), we discretize the auxiliary initial value problem
\br
& & \phi_t = \frac{1}{2} \Delta \phi -  \lambda b_1(y)b_2(t) \phi, \no \\
& &  \phi(y,0) \equiv 1
\er
using the Crank-Nicholson scheme:
\be
\frac{\phi^{n+1}_m - \phi^n_m}{\Delta t} = \frac{1}{2}(D^2 \phi^{n+1}_m + D^2 \phi^n_m) + \frac{1}{2} (F^{n+1}_m \phi^{n+1}_m + F^n_m \phi^n_m), \label{e44}
\ee
where $D^2 \phi^n_m$ denotes the standard second order discretization of the Laplacian (or $\frac{1}{2} \Delta_y$) centered at the discrete point $(y_m,t_n)$. The term $F^n_m$ corresponds to the reaction term $\lambda b_1(y)b_2(t)$ evaluated at discrete points $(y_m,t_n)$. This scheme is implicit, second-order in both time and space. In all simulations we use $b_1(y) = \delta \sin(6 \pi y)$, where $\delta >0$ is a scaling parameter, and we compute on the domain $y \in [0,1]$ with discrete grid spacing $\Delta y = 0.01$.  We use an adaptive time step, since the implicit treatment of the reaction term requires that $\Delta t < 2/(F^{n+1})$. To generate realizations of $b_2(t)$ we integrate the It\^o equation (\ref{e43}) using an implicit order 2.0 strong Taylor scheme (see \cite{KoPl}) with a discrete spacing $\Delta t_{b_2} \leq 0.1 (\Delta t)$, where $\Delta t$ is the adapting time step for the PDE evolution.

\begin{figure}[tb]
\centerline{\epsfig{file=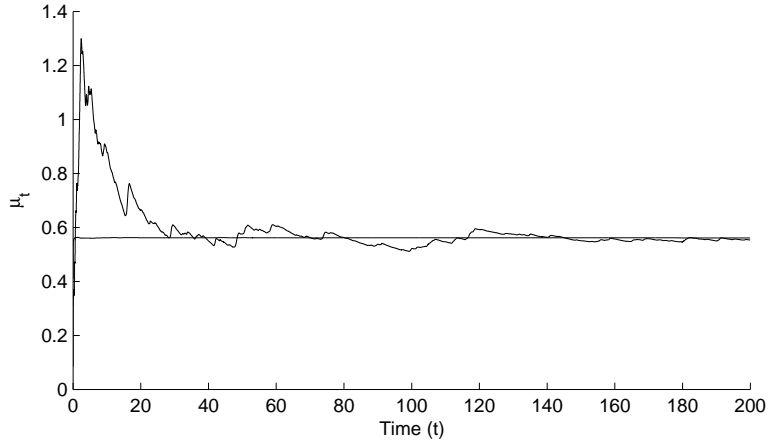,width=330pt}}
\caption{One realization of $\rho_t = \frac{1}{t}\log(\norm{\phi(y,t)}_1) = \mu_t - \lambda^2/2 - f'(0)$. The nearly flat curve shows the sample mean, $N=40,000$ realizations. }
\label{fig1}
\end{figure}

\begin{figure}[tb]
\centerline{\epsfig{file=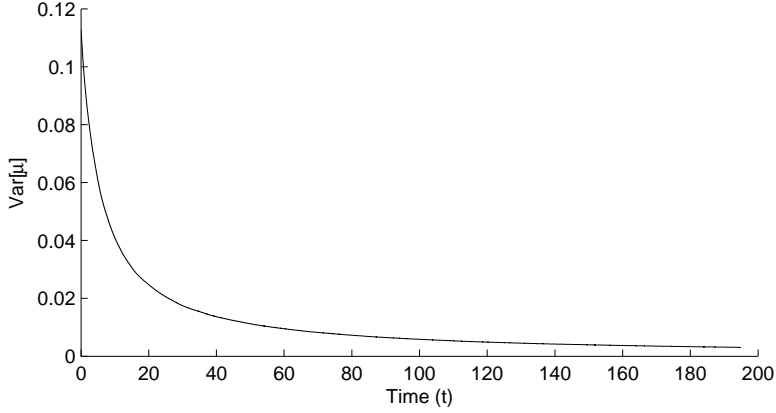,width=330pt}}
\caption{Variance of $\lambda^2/2 + f'(0) + \frac{1}{t}\log(\norm{\phi(y,t)}_1)$, $N=40,000$ realizations.}
\label{fig3}
\end{figure}

\begin{figure}[tb]
\centerline{\epsfig{file=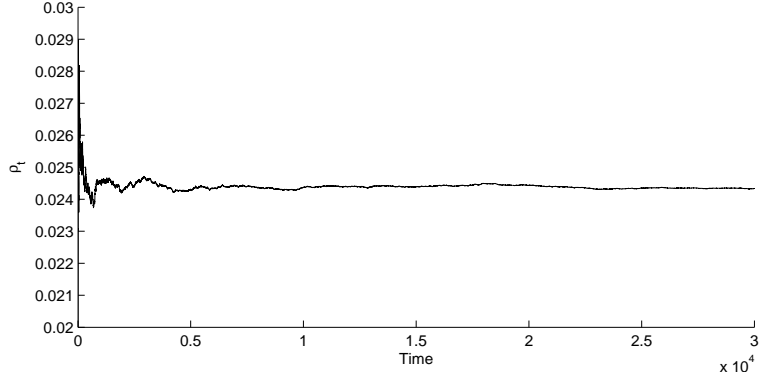,width=330pt}}
\caption{Convergence of $\rho_t = \frac{1}{t}\log(\norm{\phi(y,t)}_1) = \mu_t - \lambda^2/2 - f'(0)$. $T_f = 30,000$.}
\label{fig2}
\end{figure}

\begin{figure}[tb]
\centerline{\epsfig{file=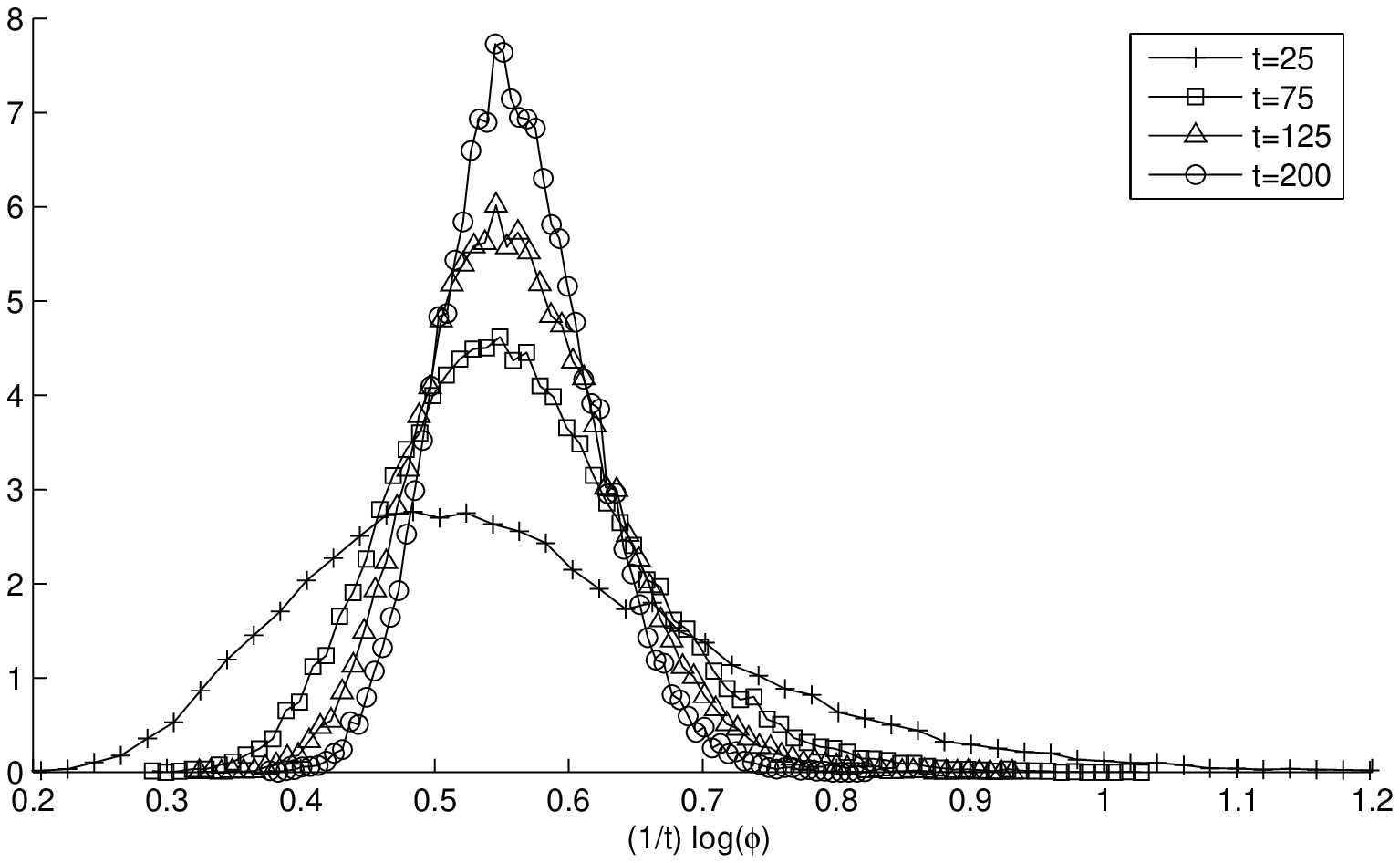,width=330pt}}
\caption{Distribution of $\frac{1}{t}\log(\norm{\phi(y,t)}_1)$ at different times.}
\label{fig4}
\end{figure}

To approximate $\mu(\lambda)$, we iterate (\ref{e44}) for a very long time $t = T_f$ and approximate
\be
\mu(\lambda) \approx \mu_t(\lambda)  =  \lambda^2/2 + f'(0) + \frac{1}{t} \log(\norm{\phi}_1). \label{e46}
\ee
By the results of the previous sections, $\mu_t$ converges to $\mu$ almost surely, so we need only generate one realization of the process. Alternatively, we could generate $N$ realizations of $b_2(t)$, evolve (\ref{e44}), and approximate
\br
\mu(\lambda) \approx E_Q[\mu_t(\lambda,\hat \omega)] & = & \lambda^2/2 + f'(0) + E_Q[\frac{1}{t} \log(\norm{\phi}_1)] \no \\
& = & \lambda^2/2 + f'(0) + \frac{1}{N} \sum_{i=1}^N \frac{1}{t} \log(\norm{\phi_i}_1).  \label{e95}
\er
In practice, we observe that $T_f$ can be chosen much smaller when using (\ref{e95}) instead of (\ref{e46}), since the mean converges must faster than an individual sample. In Figure \ref{fig1}, we show one realization of the approximation $\mu_t(\lambda)$ compared with the ensemble mean (\ref{e95}). After a very short time, the mean shows relatively little fluctuation compared to the individual realization $\mu_t(\lambda,\hat \omega)$. The variance of $\mu_t(\lambda)$, shown in Figure \ref{fig3}, decays like $O(1/t)$. However, the need to compute a large number of realizations in (\ref{e95}) makes this approach no less computationally expensive than evolving only one sample for a very large time. So, we generally use (\ref{e46}) to approximate $\mu(\lambda)$. Typically we use a final time of $T_f = 30,000$. Figure \ref{fig2} shows the convergence of $\mu$ computed by (\ref{e46}).  In Figure \ref{fig4} we show the distribution of $\mu_t$ at different points in time. For this simulation we used $N = 40,000$ realizations.

\subsection{Computation of $c^*$}
Computing the speed $c^*$ is quite simple: we use the above method to evaluate $\mu(\lambda)$ at different points on the curve $\lambda \mapsto \frac{\mu(\lambda)}{\lambda}$ and then minimize in $\lambda$.  Direct simulation of $c^*$ would require evolution of (\ref{e0}), a nonlinear, time-dependent PDE in two space dimensions. Finding $c^*$ by computing the curve $\mu(\lambda)/\lambda$, however, reduces the PDE computation to one dimension. The trade-off is that we now have a minimization problem in $\lambda$, but this is easily accomplished with a standard algorithm \cite{Min}.  From the representation (\ref{e64}), one can see that the curve $\mu(\lambda)/\lambda$ has the following properties:
\begin{lem}\label{lem:uniquemin}
The infimum of the curve $\frac{\mu(\lambda)}{\lambda}$ over $(0,\infty)$ is achieved at a unique point $\lambda^* \in (0,\lambda_0]$ where $\lambda_0 = \sqrt{f'(0)/2}$. Moreover, there are no other local minima.
\end{lem}
{\bf Proof:}
This follows from the fact that $\mu(\lambda) = \lambda^2/2 + f'(0) + \rho(\lambda)$ with $\rho$ being convex in $\lambda$ and $\rho(0) = 0$ (see discussion leading to (\ref{e97})). The point $\lambda_0$ is the value of $\lambda$ where the infimum of the curve $\lambda/2 + f'(0)/\lambda$ is attained.
\qed

Next we consider the scaling $b(y,t) \mapsto \delta b(y,t)$ and the resulting enhancement of the corresponding speed $c^* = c^*(\delta)$. It is known \cite{nolen:ekt04} that if $b(y,t)$ is periodic in both space and time that $c^*(\delta) = c^*(0) + O(\delta^2)$ for $\delta$ small and $c^*(\delta) = c^*(0) + O(\delta)$ for $\delta$ as $\delta \to \infty$. The following proposition gives analytical upper bounds consistent with this asymptotic behavior.

\begin{prop}[Bounds on $c^*$] \label{prop:ubounds}
For all $\delta \geq 0$, $c^*(\delta)$ satisfies the bounds
\bi
\item[(i)] $c^*(\delta) \geq c^*(0)$.
\item[(ii)] $c^*(\delta) \leq c^*(0) + \delta \norm{b_1}_\infty E_Q[\abs{b_2}]$.
\item[(iii)] $c^*(\delta) \leq c^*(0) \sqrt{1 + \delta^2 p_1}$.
\item[(iv)] $c^*(\delta) = c^*(0)\;\;\text{if}\;\;b(y,t) = b(t)$.
\ei
From (iii), we also have
\be
c^*(\delta) \leq c^*(0)(1 + \frac{\delta^2 p_1}{2}) + O(\delta^3) \no
\ee
when $\delta$ is small.
\end{prop}
{\bf Proof:} The first bound follows from (\ref{e97}) and the formula
\be
c^*(\delta) = \inf_{\lambda > 0} \frac{\mu(\lambda)}{\lambda} \geq \inf_{\lambda > 0} \frac{\lambda}{2} + \frac{f'(0)}{\lambda} = c^*(0). \no
\ee
The function $\psi = \log(\phi)$ satisfies
\br
\psi_t &=& \frac{1}{2} \Delta \psi + \frac{1}{2}\abs{\nabla_y \psi}^2 - \lambda b(y,t) \label{e87} \\
\psi(y,0) &\equiv & 0. \no
\er
Integrating (\ref{e87}) over $D \times [0,t]$, we have
\be
\frac{1}{t} \int_D \psi(y,t) \,dy =  \frac{1}{2t} \int_0^t \int_D  \abs{\nabla_y \psi}^2\,dy\,dt - \frac{\lambda}{t} \int_0^t \int_D b(y,t) \,dy\,dt. \label{e103}
\ee
Now let $t \to \infty$:
\br
\rho(\lambda) = \lim_{t \to \infty} \frac{1}{t} \int_D \psi(y,t) \,dy &\geq& - \lambda \lim_{t \to \infty} \frac{1}{t} \int_0^t \int_D b(y,t) \,dy\,dt \no \\
& = & - \lambda E_Q\left[ \int_D b(y,t) \,dy \right] = 0, \no
\er
almost surely with respect to $Q$.  If $b(y,t) = b(t)$, then the first integral on the right hand side of (\ref{e103}) vanishes since $\abs{\nabla_y \psi}^2 \equiv 0$. Then taking the limit as $t \to \infty$ we have equality:
\be
\rho(\lambda) = E_Q\left[ \int_D \delta b(t) \,dy \right] = 0. \no
\ee
Hence $c^*(\delta) = c^*(0)$. This proves (iv).

For the linear upper bound (ii), note that
\br
\frac{1}{t} \log E[e^{\lambda \delta \int^t_0 b_1(W(s))b_2(t-s)\,ds}] & \leq & \frac{1}{t} \log E[e^{\lambda \delta \norm{b_1}_\infty \int^t_0 \abs{b_2(t-s)}\,ds}] \no \\
& = & \frac{\abs{\lambda} \delta \norm{b_1}_\infty}{t} \int_0^t \abs{b_2(s)} \,ds. \no
\er
As $t \to \infty$, this last term converges almost surely to $\abs{\lambda} \delta \norm{b_1} E_Q[\abs{b_2}]$. Therefore, $c^*(\delta)$ always satisfies the linear upper bound
\br
c^*(\delta) = \inf_{\lambda > 0} \frac{\mu(\lambda)}{\lambda} & \leq & \inf_{\lambda > 0} \frac{\lambda}{2} + \frac{f'(0)}{\lambda} + \delta \norm{b_1} E_Q[\abs{b_2}] \no \\
& = & c^*(0) + \delta \norm{b_1} E_Q[\abs{b_2}]. \label{e90}
\er
Finally, for the quadratic upper bound, observe that under the scaling $b \mapsto \lambda \delta b$, the constant $p_1$ defined in A5 can be replaced by $p_1 \mapsto \lambda^2 \delta^2 p_1$. Then by (\ref{e85}) and (\ref{e84}),
\be
\rho(\lambda) \leq \sqrt{2} \lambda^2 \delta^2 p_1 \no
\ee
and
\br
c^*(\delta) = \inf_{\lambda > 0} \frac{\mu(\lambda)}{\lambda} & \leq & \inf_{\lambda > 0} \frac{\lambda}{2} + \frac{f'(0)}{\lambda} + \frac{\lambda^2 \delta^2 p_1}{2} = c^*(0) \no \\
& = & 2 \sqrt{(1 + \delta^2 p_1)f'(0)/2} \no \\
& = & c^*(0) \sqrt{(1 + \delta^2 p_1)} \label{e92} \\
& = & c^*(0)(1 + \frac{\delta^2 p_1}{2}) + O(\delta^3). \no
\er
\qed

If we make some additional restriction on the form of $b(y,t)$ we also have a linear lower bound on the growth of $c^*(\delta)$ as $\delta \to \infty$.
\begin{prop}[Linear growth of $c^*$]
Let $b(y,t)$ have the form
\be
b(y,t) = \sum_{j=1}^{N} b^j_1(y) b^j_2(t) \no
\ee
where $b^j_1(y)$ are Lipschitz continuous and periodic in y, and $b^j_2(t)$ are stationary centered Gaussian fields such that the Assumptions A1-A6 are satisfied for $b(y,t)$. Then the constant $\bar C \in [0 ,+\infty)$ defined by
\be
\liminf_{\delta \to \infty} \frac{c^*(\delta)}{\delta} = \bar C
\ee
is equal to zero if and only if $b(y,t) \equiv b(t)$.
\end{prop}
{\bf Proof:}
The fact that $\bar C \in [0,+\infty)$ follows from Proposition \ref{prop:ubounds}. Also, if $b(y,t) \equiv b(t)$ then $\bar C = 0$ since $c^*(\delta) = c^*(0)$ for all $\delta > 0$. By Lemma \ref{lem:uniquemin} there is a unique $\lambda = \lambda_\delta \in (0,\lambda_0]$ such that
\be
c^*(\delta) = \inf_{\lambda > 0} \frac{\mu(\lambda)}{\lambda} = \frac{\mu(\lambda_\delta)}{\lambda_\delta}. \no
\ee
Let $\delta_j \to \infty$ as $j \to \infty$ and suppose that $\limsup_{j \to \infty} (\lambda_{\delta_j} \delta_j) \leq M$. This implies that
\be
\liminf_{j \to \infty} \frac{c^*(\delta_j)}{\delta_j} = \liminf_{j \to \infty} \frac{\mu(\lambda_{\delta_j})}{\lambda_{\delta_j} \delta_j} \geq \liminf_{j \to \infty} \frac{f'(0)}{\lambda_{\delta_j} \delta_j} \geq \frac{ f'(0)}{M} > 0. \no
\ee
So, in this case the result holds with $\bar C = f'(0)/M$.

Now suppose $\lambda_{\delta_j} \delta_j$ is unbounded as $j \to \infty$. By Proposition \ref{prop:rholinear}, there is a positive constant $K$ such that
\be
\rho(\lambda_{\delta_j} \delta_j) \geq K \lambda_{\delta_j} \delta_j > 0 \no
\ee
for $j$ sufficiently large. Note that Proposition \ref{prop:rholinear} treats the case of $\delta = 1$; this is why we use $\rho(\lambda_{\delta_j} \delta_j)$ instead of $\rho(\lambda_{\delta_j})$. Therefore,
\br
\liminf_{j \to \infty} \frac{c^*(\delta_j)}{\delta_j} & = & \liminf_{j \to \infty} \frac{ \lambda_{\delta_j}}{2 \delta_j} + \frac{f'(0)}{\lambda_{\delta_j} \delta_j} + \frac{\rho(\lambda_{\delta_j} \delta_j) }{ \lambda_{\delta_j} \delta_j} \geq K > 0 \no
\er
since $\lambda_{\delta_j} \in (0,\lambda^*]$ and $\lambda_{\delta_j} \delta_j \to \infty$. Hence $\bar C \geq K > 0$.
\qed

\begin{figure}[tb]
\centerline{\epsfig{file=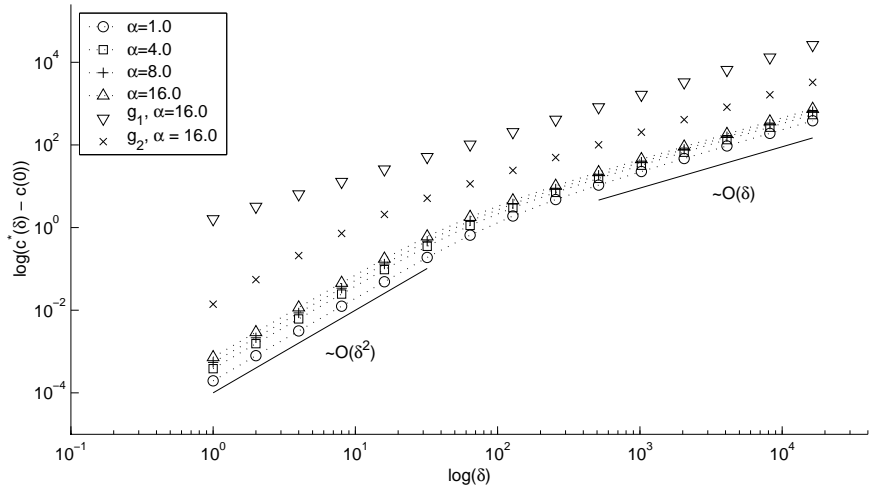,width=330pt}}
\caption{Log-log plot of front speed $c^*(\delta)$ versus $\delta$. For comparison, the solid lines on the bottom have slope $p=2.0$ and $p=1.0$. The data sets $g_1$ and $g_2$ represent the upper bounds (\ref{e90}) and (\ref{e92}), respectively, contained in Proposition \ref{prop:ubounds}. Both bounds were computed for $\alpha = 16.0$.}
\label{fig10}
\end{figure}

The results of our numerical computations suggest that $c^*(\delta)$ is a monotone increasing function of $\delta$, and they confirm both the quadratic and linear growth of $c^*(\delta)$ for small and large $\delta$, respectively.  In Figure \ref{fig10}, we plot $c^*(\delta) - c^*(0)$ on a log-log scale for a few values of the covariance parameter $\alpha$. We observe a transition form quadratic scaling in the small $\delta$ regime to linear scaling in the large $\delta$ regime. We also plot the upper bounds $g_1$ and $g_2$ given by (ii) and (iii) of Proposition \ref{prop:ubounds} using $\alpha = 16.0$.  We observe that these bounds are at least an order of magnitude greater than the numerically computed enhancement.

Next, by varying the parameter $\alpha$ in (\ref{e91}), we consider the effect of temporal correlations on the enhancement of the speed $c^*$. Because of the simple structure of the field $b_1(y)b_2(t)$ we can compute the constant $p_1$ appearing in Assumption A5 and Proposition \ref{prop:ubounds}. Using (\ref{e91}),
\br
\sup_{y_1,y_2} \Gamma(y_1,y_2,0,r) = \norm{b_1}_\infty^2 E_Q[b_2(0)b_2(r)] = \sqrt{\alpha} \norm{b_1}_\infty^2 e^{-\alpha r}. \no
\er
Hence $p_1 = \alpha^{-1/2} \norm{b_1}_\infty^2 $. Therefore, from (iii) of Proposition \ref{prop:ubounds},
\be
c^* \leq c^*(0) \sqrt{(1 + \delta^2 \alpha^{-1/2} \norm{b_1}_\infty^2 )}  \label{e93}
\ee
So, as $\alpha \to \infty$, $c^* \to c^*(0)$. This limit corresponds to the correlation length $1/\alpha$ becoming very small and is consistent with the case of periodic time dependence: faster temporal oscillation of the shear tends to decrease the enhancement of the front speed \cite{Kh, NX1, nolen:ekt04}.  As $\alpha \to 0$, the bound (\ref{e93}) blows up. From (ii) of Proposition \ref{prop:ubounds}, however, we also have
\br
c^* & \leq & c^*(0) + \delta \norm{b_1} E_Q[\abs{b_2}] \no \\
 & = &  c^*(0) + \delta \alpha^{1/4} \norm{b_1} E_Q[\abs{Z}]  \label{e94}
\er
where $Z$ is a normally distributed random variable. This implies that $c^* \to c^*(0)$ as $\alpha \to 0$, as well. So for $\alpha \in (0, \infty)$, there must be some optimal correlation length $1/\alpha$ so that enhancement is maximal.

\begin{figure}[tb]
\centerline{\epsfig{file=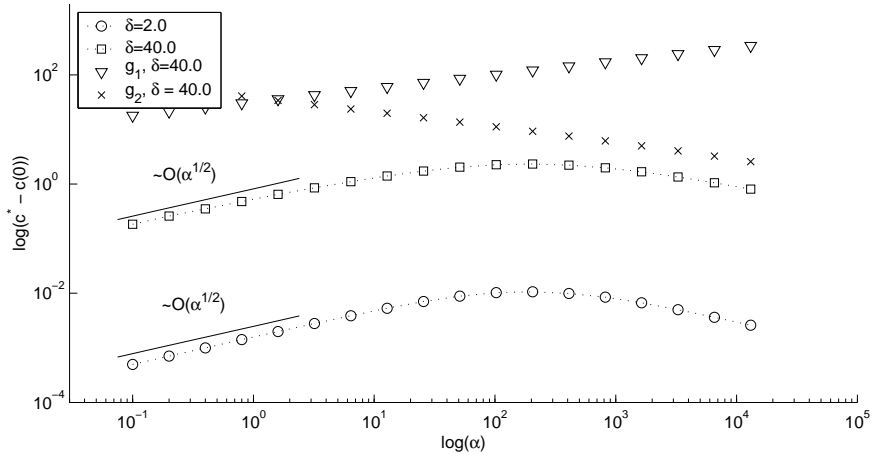,width=330pt}}
\caption{Log-log plot of $c^*$ vs. covariance parameter $\alpha$ ($\delta = 2.0, 40.0$). For comparison, the solid lines have slope $p=1/2$. The data sets $g_1$ and $g_2$ represent the upper bounds (\ref{e90}) and (\ref{e92}), respectively, contained in Proposition \ref{prop:ubounds}. Both bounds were computed for $\delta = 40.0$.}
\label{fig11}
\end{figure}

For $\alpha \in [1/10, 2^{16}/10]$, we computed the expected speed $c^*(\delta)$ for fixed amplitudes $\delta = 2.0$ and $\delta = 40$. Note that for each $\alpha$, we must
choose the initial points $b_2(0)$ to have variance $E[b_2(0)^2]= \sqrt{\alpha}$ so that the
process remains stationary for each $\alpha$.  Also, as $\alpha$ becomes large, we adjust the PDE time step so that $\Delta t \leq 0.5/\alpha$, in addition to other restrictions already mentioned. Otherwise, the numerical method cannot resolve the fast oscillations of the shear process, and the speeds diverge as $\alpha$ grows.

\begin{figure}[tb]
\centerline{\epsfig{file=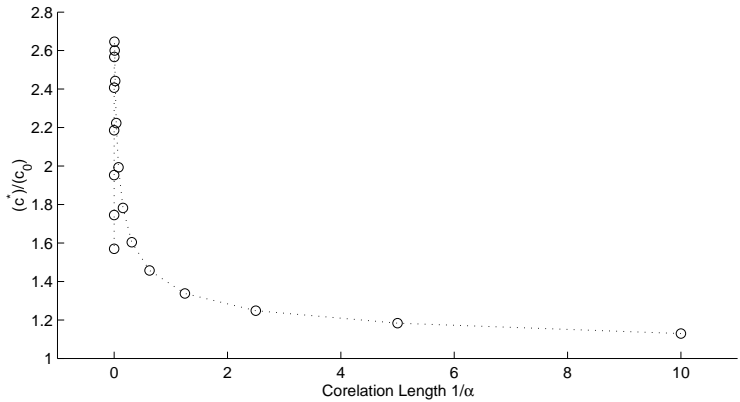,width=280pt}}
\caption{Dependence of $c^*$ on correlation length $1/\alpha$ ($\delta = 40.0$).}
\label{fig12}
\end{figure}

In Figure \ref{fig11}, we plot $c^*$ versus $\alpha$ on a log-log scale. We also plot the
upper bounds $g_1$ and $g_2$ given by (\ref{e94}) and (\ref{e93}) for $\delta=40$.
In general, these bounds are rather coarse, an order of magnitude larger than the numerically computed enhancement, but as $\alpha \to \infty$, the bound $g_2$ lies relatively close to the data. As $\alpha \to 0$, $g_1$ is the better bound, predicting that $c^* \to c^*(0)$, although the scaling of the bound is different from the scaling observed in the data. By computing the slope of a best-fit line through the points on the log-log plot, we find that for $\delta=2.0$, $c^*$ scales according to $c^* = c^*(0) + O(\alpha^{0.50})$ and for $\delta = 40$, $c^* = c^*(0) + O(\alpha^{0.44})$. So in both cases, $\abs{c^* - c^*(0)} \to 0$ much faster than the $O(\alpha^{0.25})$ convergence predicted by the bound $g_1$ in (\ref{e94}). This scaling behavior can be understood by analogy with the periodic case. Note that $\delta b_1(y) b_2(t)$ = $\alpha^{1/4} \delta b_1(y) \hat b_2(t)$ where $\hat b_2$ has unit variance and correlation length $1/\alpha$. If $b(y,t) = b_1(y)\hat b_2(t)$ and $\hat b_2$ is periodic with very long wavelength, then the enhancement is approximately equal to the enhancement caused by $b_1$ only (a steady shear). The very slow oscillations in the shear field do not significantly slow the front. In the random case, when $\alpha$ is small the correlation length is very large.  So by analogy, we expect that in the random case, when $\alpha$ is very small, $c^*$ will behave as if the shear were just $ \alpha^{1/4} \delta b_1$. In the small amplitude regime, $c^*$ scales quadratically with amplitude. Hence $c^* \approx c^*(0) + O((\alpha^{1/4} \delta)^2) = c^*(0) + O(\alpha^{1/2} \delta^2)$, which is consistent with our numerical computations for $\alpha$ small. Figure \ref{fig12} shows the data for $\delta = 40.0$ in terms of correlation length $1/\alpha$.

\section{Conclusions}
We have considered the propagation of KPP reaction fronts in
temporally random shear flows inside an infinite cylinder. We showed that,
under assumptions A1-A6 on a Gaussian shear field,
the fronts speeds obey a variational formula that extends the known
variational formula in the case of periodic media.
We performed analysis and computation of front speeds based on
the variational principle.
When the shear field is Ornstein-Uhlenbeck in time, we
numerically demonstrated the quadratic and linear speed growth laws in the shear root mean square
amplitude, and speed dependence on the shear temporal correlation length.
We also derived basic bounds on the front speeds and compared them with the
computed speeds. Developing methods to improve these bounds will be left as a future work.

\section{Acknowledgements}
The work was partially supported by NSF grants ITR-0219004, SCREMS-0322962, DMS-0506766.
J. X. would like to thank Prof. M. Cranston for helpful communications.
J. N. is grateful for support through a VIGRE graduate fellowship at UT Austin.

\section{Appendix A: Estimates on $X^{z,t}(s)$}
\setcounter{equation}{0}
In this section we derive some technical estimates on the family of processes $X^{z,t}(s)$ that follow from our assumptions on the field $b$ and the Borell inequality for Gaussian fields. In particular, we prove Lemma \ref{corol:taubound} which is needed for the bounds (\ref{e16}) and (\ref{e45}). Let use first note that by changing variables $r = s-t$, $v = s+t$, it is easy to see from Assumption A5 that
\br
\int_0^T \int_0^T \sup_{y_1,y_2} \Gamma(y_1,y_2,s,t) \,ds\,dt &\leq& 2 \int_0^{\sqrt{2}T} \int_0^{T/\sqrt{2}} \hat \Gamma(r) \,dr\,dv \no \\
&\leq & 2 \sqrt{2} p_1 T \label{e88}
\er
and for $H \in [0, T]$,
\br
\int_0^T \int_H^T \sup_{y_1,y_2} \Gamma(y_1,y_2,s,t) \,ds\,dt &\leq & \sqrt{2} \abs{T-H} \int_0^{T/\sqrt{2}} \hat \Gamma(r) \,dr \no \\
& \leq & \sqrt{2} \abs{T-H}p_1. \label{e89}
\er
Let $\rho(s) \in C([0,+\infty), \mathbb{R}^{n-1})$ with $\rho(0) = 0$ be fixed. For $y \in D$, define $\rho_y(s) = y + \rho(s)$. For fixed $t >0 $, the integral
\be
f(y,s) = \int_0^s b(\rho_y(\tau),t-\tau)\,d\tau \label{e15}
\ee
is a Gaussian random field over $M = D \times [0,t]$, with respect to the measure $\hat Q$. The Borell inequality for Gaussian fields states that if $\norm{f} = \sup_{(x,s) \in M} f(x,s) $ is almost surely finite, then $E[\norm{f}] < \infty$ and for any $u>0$,
\be
Q\left( \norm{f} - E[\norm{f}] > u \right) \leq e^{-\frac{u^2}{2 \sigma^2_t}} \label{e14}
\ee
where $\sigma^2_t = \sup_{(x,s) \in M} E_Q[f^2]$ (see \cite{Adler2}). By (\ref{e88}), $\sigma^2_t \leq 2\sqrt{2}p_1 t$. Using inequality (\ref{e14}), we can control deviations of $\norm{f}$, if we bound the growth of $E[\norm{f}]$.
\begin{lem} \label{lem:expectbound}
There is a finite constant $C >0$ such that
\be
E[\norm{f}] \leq Ct^{1/2} .
\ee

\end{lem}
{\bf Proof:} The expectation $E[\norm{f}]$ can be bounded by the metric entropy relation \cite{Adler2}
\be
E[\norm{f}] \leq C \int_0^\delta \left(\log N(\epsilon) \right)^{1/2} \,d\epsilon \no
\ee
where $\delta = diam(M)/2$ in the metric
\be
d((x,s), (y,z)) = E\left[ (f(x,s) - f(y,z))^2 \right]^{1/2} \no
\ee
and $N(\epsilon)$ is the minimum number of $\epsilon$ balls required to cover $M$. Using (\ref{e88}), (\ref{e89}), and Assumption A6, a straightforward computation shows that
\be
E\left[ (f(x,s) - f(y,s))^2 \right] \leq C \abs{x-y} t \no
\ee
and
\be
E\left[ (f(y,s) - f(y,z))^2 \right] \leq C \abs{s-z} \no
\ee
for some finite constant $C$, independent of $\rho$. Therefore,
\be
d((x,s), (y,z)) \leq C_1 \left( \abs{s-z})^{1/2} + C_2 ( \abs{x-y} t \right) ^{1/2}, \no
\ee
and there is a constant $C_3$ independent of $t$ and $\epsilon$ such that $ d((x,s), (y,z)) \leq \epsilon $
whenever $\abs{s-z} \leq C_3 \epsilon^2 $ and $\abs{x-y} \leq \frac{C_4 \epsilon^2}{t}$. For $\epsilon \in (0,diam(M)/2]$, we have the bound
\be
N(\epsilon ) \leq  \max( C_5\frac{t^2}{\epsilon^4},1) \no
\ee
and
\br
E[\norm{f}] &\leq& C\int_0^{C_5^{-1/4} t^{-1/2}} (\log ( C_5 \frac{t^2}{\epsilon^4}) )^{1/2} \,d\epsilon \no \\
&  = & C_6 t^{1/2}  \int_0^1 (\log(\frac{1}{\epsilon^4}))^{1/2},d\epsilon \leq C_7 t^{1/2}.  \label{e38}
\er
\qed

Note that the constants depend on the assumed properties of the process $b$ and the size of the domain $D$, but not on the particular function $\rho(s)$. If $u \geq 2 C_7 t^{1/2}$, then by (\ref{e14}),
\be
Q \left( \norm{f} > u \right) \leq e^{-(u - E \norm{f})^2 / 2 \sigma^2_t} \leq e^{-u^2 / 8 \sigma^2_t} \leq e^{-u^2 / 8 p_1 t}. \no
\ee
It now follows that
\begin{lem}\label{lem:borelbound}
For any $\eta > 0$ and for $t \geq t_0 = t_0(\eta) = (2 C_7 / \eta)^2$,
\be
Q\left( \sup_{y \in D,\; s \in [0,t]} \int_0^s b(y + \rho(\tau),t-\tau)\,d\tau  > \eta t \right) \leq e^{- \eta^2 t / 8 p_1}. \no
\ee
for any $\rho \in C([0,\infty),R^n)$, $\rho(0) = 0$.
\end{lem}
In applying this lemma, the continuous function $\rho$ will be a realization of the Wiener process $W^y_2(s)$.

\begin{lem} \label{lem:taubound}
For $\eta > 0$, $z \in R^n$, define the Markov time
\be
\tau_{\eta,z} (t) = \min \{ s \geq 0 | \; \abs{X^{z,t}(s) - x} \geq \eta t \}. \no
\ee
with $\tau_{\eta,z} (t) = +\infty$ if the set on the right is empty. Then there are constants $K_1, K_2$ such that
\be
Q\left( P \left( \inf_{z \in R^n} \tau_{\eta,z} (t) \leq t \right) > e^{- K_2 \eta^2 t /2} \right) \leq K_1 e^{-K_2 \eta^2 t/ 2} \no
\ee
for all $t > 0$.
\end{lem}

{\bf Proof:}
Note that for the Wiener process $W_1(s)$ with $W_1(0) = 0$
\br
P \left( \sup_{s \in [0,t] } \abs{W_1(s)} \geq \eta t \right) & \leq & 2 \sqrt{\frac{2}{\pi}} \int_{\eta \sqrt{t}}^{\infty} e^{-x^2 /2}\;dx \no \\
& \leq & K_1 e^{-\eta^2 t/2} . \label{e33}
\er
The point of the lemma is that at large times, and almost surely with respect to $Q$, the process $X^{z,t}(s)$ still behaves like a Wiener process, even though the drift term $b(y,t)$ is not uniformly bounded in $t$. By definition of $\tau_{\eta,z}(t)$,
\br
P \left( \inf_{z \in R^n} \tau_{\eta,z} (t) \leq t \right) &=& P\left( \sup_{s\in [0,t], z \in R^n} \abs{f_t(y,s) + W_1(s)} \geq \eta t  \right). \no
\er
Using Tchebyshev's inequality, (\ref{e33}), and Lemma \ref{lem:borelbound} we see that for any $\eta > 0$, $\alpha > 0$:
\br
Q\left( P \left( \inf_{z \in R^n} \tau_{\eta,z} (t) \leq t \right) > \alpha \right) & \leq & \alpha^{-1} E_Q P\left( \sup_{s\in [0,t], z \in R^n} \abs{f_t(y,s) + W_1(s)} \geq \eta t  \right) \no \\
& =& \alpha^{-1} E_P Q \left( \sup_{s\in [0,t], z \in R^n} \abs{f_t(y,s) + W_1(s)} \geq \eta t  \right) \no \\
& \leq &  \alpha^{-1} E_P Q \left( \sup_{s\in [0,t], z \in R^n} \abs{f_t(y,s)} \geq \eta t/2  \right) \no \\
& & +  \alpha^{-1} P \left( \sup_{s\in [0,t], z \in R^n} \abs{ W_1(s)} \geq \eta t/2  \right) \no \\
& \leq & \alpha^{-1} (2 e^{-\eta^2 t/32 p_1} + K_1 e^{- \eta^2 t/8} ) \leq   \alpha^{-1} K_1 e^{- K_2 \eta^2 t } \no
\er
for $t$ sufficiently large, for some constants $K_1, K_2 > 0$. The result now follows from a choice of $\alpha = e^{- K_2 \eta^2 t /2}$.
\qed
\begin{cor} %\label{corol:taubound}
There are constants $K_1, K_2 > 0$ such that, except on a set of $ Q$-measure zero,
\be
\sup_{ z \in R^n} P\left( \tau_{\eta,z} (t) \leq t \right) \leq K_1 e^{- K_2 \eta^2 t}
\ee
for $t$ sufficiently large depending on $\hat \omega$ and $\eta$.
\end{cor}
{\bf Proof:}
Note that this is Lemma \ref{corol:taubound} in the case that $\kappa = 1$. Lemma \ref{lem:taubound} and the Borel-Cantelli lemma imply that outside a set of $Q$-measure zero
\be
P\left( \inf_{z \in R^n} \tau_{\eta,z} (k) \leq k \right) \leq e^{- K_2 \eta^2 k /2} \label{e34}
\ee
if $k \in \mathbb{Z}$ is sufficiently large. Now we want to extend this to all real $t$ sufficiently large. Let $t \in [k, k+1]$, $t = k + \tau$, $\tau \in [0,1]$.
\br
& & \sup_{z \in R^n} P\left( \sup_{s \in [0,t]} \abs{X^{z,t}(s) - x_0} \geq t \eta  \right) \no \\
& \leq &   \sup_{z \in R^n}P\left( \sup_{s \in [0,\tau]} \abs{X^{z,t}(s) - x_0} \geq t \eta/2 \right) +  \no \\
& & \quad \quad + \sup_{z \in R^n}P \left( \sup_{s \in [\tau,t]} \abs{X^{z,t}(s) - X^{z,t}(\tau)} \geq t \eta/2 \right)\no
\er
By the Markov property, this is bounded by
\br
&  \leq & \sup_{z \in R^n}P\left( \sup_{s \in [0,\tau]} \abs{X^{z,t}(s) - x_0} \geq t \eta/2 \right) +  \no \\
& & \quad \quad + \sup_{\bar z \in R^n}P\left( \sup_{s \in [0,k]} \abs{X^{\bar z,k}(s) - \bar x_0} \geq t \eta/2   \right) \leq \no \\
& \leq &  P\left( \sup_{z \in R^n, t \in [k,k+1], s \in [0,1]} \abs{X^{z,t}(s) - x_0} \geq k \eta/2 \right) +  \no \\
& & \quad \quad + P\left(  \sup_{\bar z \in R^n, s \in [0,k]} \abs{X^{\bar z,k}(s) - \bar x_0} \geq k \eta/2   \right) \label{e36}
\er
By (\ref{e34}), the second term on the right side of (\ref{e36}) is bounded (Q-a.s.) by
\br
P\left(  \sup_{\bar z \in R^n, s \in [0,k]} \abs{X^{\bar z,k}(s) - \bar x_0} \geq k \eta/2   \right) \leq e^{- K_3 \eta^2 k } \label{e37}
\er
for $k \in \mathbb{Z}$ sufficiently large. To bound the other term in (\ref{e36}), it suffices to show that
\be
 P\left( \sup_{y \in D,  r \in [0,1], s \in [0,1]} \abs{f_k(y,s,r)} \geq k \eta/2 \right) \leq e^{- K_4 \eta^2 k } \label{e35}
\ee
for $k \in \mathbb{Z}$ sufficiently large, where
\be
f_k(y,s,r) = \int_0^s b(W^y_2(\tau),r+k-\tau)\,d\tau . \no
\ee
Note that $f_k(y,s,r)$ is a centered Gaussian field over $D \times [0,1] \times [0,1]$, and its distribution is invariant with respect to $k > 0$, due to the stationarity of $b(y,t)$. For any fixed path $W^y_2(\omega)$, the Borell inequality implies that for $k$ sufficiently large
\be
Q\left( \sup_{y \in D,  r \in [0,1], s \in [0,1]}  \abs{f_k(y,s,r)} \geq k \eta/2 \right) \leq K_5 e^{- K_6 \eta^2 k^2} \no
\ee
for some constants $K_5, K_6 > 0$, independent of $k$ and the realization $W^y_2(\omega)$. Therefore, proceeding as in the proof of Lemma \ref{lem:borelbound}, we see that
\be
Q\left( P( \sup_{y \in D,  r \in [0,1], s \in [0,1]} \abs{f_k(y,s,r)} \geq k \eta/2 ) \geq e^{-K_6 \eta^2 k^2/2} \right) \leq K_7 e^{- K_6 \eta^2 k^2 /2} . \no
\ee
Now (\ref{e35}) follows from the Borel-Cantelli lemma. We complete the proof by combining (\ref{e37}) and (\ref{e35}).
 \qed

{\bf Proof of Lemma \ref{corol:taubound}:}
We have just proved Lemma \ref{corol:taubound} in the special case that $\kappa = 1$. For $\kappa < 1$,  modify the preceding bounds for the field
\be
f(y,s) = \int_0^s b(\rho_y(\tau),t-\tau)\,d\tau \no
\ee
considered over $M_{\kappa} = D \times [0,\kappa t]$. Now we have $\sigma^2_t = \sup_{(x,s) \in M} E_Q[f^2] \leq p_1 \kappa t$, so we find that $E[\norm{f}] \leq C \sqrt{\kappa t}$ for some constant $C>0$. Then, just as in Lemma \ref{lem:taubound}, we have
\be
Q\left( P \left( \sup_{s \in [0,\kappa t], z \in R^n} \abs{X^{z,t}(s) - x} \geq \eta t  \right) > e^{- K_2 \eta^2 t /2\kappa } \right) \leq K_1 e^{-K_2 \eta^2 t/ 2 \kappa}. \no
\ee
The rest follows as in the preceding corollary.
\qed

\section{Appendix B: Lyapunov Exponent $\mu(\lambda)$}
\setcounter{equation}{0}
In this section we prove Lemma \ref{lem:muexist} assuming that A1-A6 hold for the random field $b(y,t)$.  We study the limit
\br
\mu(\lambda,z)  &=&  f'(0) + \lim_{t \to \infty} \frac{1}{ t} \log E\left[ e^{  - \lambda X^{z,t}(t) } \right] \no \\
& = & f'(0) + \frac{\lambda^2}{2} + \lim_{t \to \infty} \frac{1}{ t} \log E_{y} \left[ e^{ - \lambda \int_0^{t} b(W^y_2(s),t-s)\,ds } \right] \no \\
& = & f'(0) + \frac{\lambda^2}{2} + \rho(\lambda,y). \label{e96}
\er
By the Feynman-Kac formula
\br
\rho(\lambda, y) & = &  \lim_{t \to \infty} \frac{1}{t} \log \phi(y,t) \label{e64}
\er
where $\phi(y,t) > 0$ solves that auxiliary initial value problem
\br
\phi_t = \frac{1}{2} \Delta \phi -   \lambda b(y,t) \phi \label{e86} \\
\phi(y,0) \equiv 1. \no
\er
The equation (\ref{e86}) is called the parabolic Anderson problem (see \cite{CarMol} and \cite{CM}). The proof that $\rho(\lambda)$ exists (and $\mu(\lambda)$) almost surely with respect to $Q$, independent of $z$, relies on the subadditive ergodic theorem and a recent result in \cite{CM}, provided we assume the necessary decay of time correlation of the process $b(y,t)$. For simplicity of notation we ignore the dependence on $-\lambda$, $-\lambda b \mapsto b$. Following \cite{CM}, we define for any continuous path $X \in C([0,t],R)$ the exponential
\be
\xi(t,X) = e^{\int_0^t b(X_s,t-s)\,ds }.  \no
\ee
For any fixed path $X$, $\xi(t,X)$ is lognormal with mean and variance
\be
E_Q [\xi(t,X)] = e^\frac{\hat \sigma^2}{2},\;\;\;\; Var[\xi(t,X)] = e^{2 \hat \sigma^2} - e^{\hat \sigma^2} \no
\ee
where
\be
\hat \sigma^2 = \int_0^t \int_0^t \Gamma( X_s, X_r, s, r) \,ds \,dr \leq  2 \sqrt{2} p_1 t , \label{e110}
\ee
by (\ref{e88}). Note that $\hat \sigma^2$ is bounded independently of the particular path $X$.

For $0 \leq s < t$, define the random variables
\br
q_I(s,t) &=& \inf_y E_y[e^{\int_0^{t-s} b(W(\tau),t-\tau)\,d\tau}]  \no \\
q_S(s,t) &=& \sup_y E_y[e^{\int_0^{t-s} b(W(\tau),t-\tau)\,d\tau}]. \no
\er
Using the subadditive ergodic theorem, we will show that the limits
\be
\lim_{t \to \infty} \frac{1}{t} \log q_I(0,t) = \rho_I \label{e111}
\ee
and
\be
\lim_{t \to \infty} \frac{1}{t} \log q_S(0,t) = \rho_S \label{e112}
\ee
exists and are finite, almost surely with respect to $Q$. Then we will show $\rho_I = \rho_S$. By the Markov property of the Wiener process we have for and $s < z < t$:
\br
q_I(s,t) & = &\inf_y E_y\left[e^{\int_0^{t-z} b(W(\tau),t-\tau)\,d\tau} e^{\int_{t-z}^{t-s} b(W(\tau),t-\tau)\,d\tau}\right] \no \\
& = &\inf_y E_y\left[e^{\int_0^{t-z} b(W(\tau),t-\tau)\,d\tau} E[e^{\int_{t-z}^{t-s} b(W(\tau),t-\tau)\,d\tau}| W_{t-z}]\right] \no \\
& \geq &\inf_y E_y\left[e^{\int_0^{t-z} b(W(\tau),t-\tau)\,d\tau} \inf_q E_q[e^{\int_{0}^{z-s} b(W(\tau'),z-\tau')\,d\tau'}]\right] \no \\
& = & q_I(s,z)q_I(z,t). \no
\er
Therefore, $\log(q_I(s,t))$ is super-additive:
\be
\log(q_I(s,t)) \geq \log(q_I(s,z)) + \log(q_I(z,t)) \no
\ee
for any $0 \leq s < z < t$. Similarly, the function $\log(q_S(s,t))$ is subadditive. By the stationarity of $b$,
\be
\tau_z \log(q_I(s,t)) = \log(q_I(s+z,t+z)) \no
\ee
for any $z \geq 0$. Moveover, $\log(q_I(s,t))$ is integrable:
\br
E_Q[ \log(q_I(s,t)) ] \leq \log E_Q[q_I(s,t)] & \leq & \log \inf_y E_P E_Q [ e^{\int_0^{t-s} b(W(\tau),t-\tau)\,d\tau} ] \no \\
& \leq & \log e^{\frac{\sigma^2}{2}} \leq \sqrt{2} p_1 \abs{t-s}, \label{e85}
\er
and
\br
E_Q\left[ \log(q_I(s,t)) \right] & \geq & E_Q\left[ \log( E_P e^{\inf_y \int_0^{t-s} b(W^y(\tau),t-\tau)\,d\tau} ) \right] \no \\
& \geq & E_P E_Q \left[  \inf_y \int_0^{t-s} b(W^y(\tau),t-\tau)\,d\tau ) \right].
\er
This last term is finite, by the Borell inequality. Note also that
\br
\sup_{s,t \in [0,1]} \abs{\log q_I(s,t)} \leq  \sup_{t \in [0,1], y \in D} \abs{b(y,t)}, \no
\er
and the latter is integrable with respect to $Q$. It now follows from the subadditive ergodic theorem (Theorem 2.5 of \cite{AK}) and the continuity of $q(0,t)$ with respect to $t$ that the limit (\ref{e111}) exists almost surely and is finite:
\be
\lim_{t \to \infty} \frac{1}{t} \log \inf_{y} \phi(y,t) = \lim_{t \to \infty} \frac{\log(q_I(0,t))}{t} = \sup_t \frac{E_Q \left[\log(q_I(0,t))\right]}{t} = \rho_I \label{e39}
\ee
Also, by (\ref{e85}), $\rho_I \leq \sqrt{2} p_1$.

In order to apply the ergodic theorem to $(1/t)\log(q_S(0,t)$, we need to show that $(1/t)\log(q_S(0,t)$ is integrable. By Jensen's inequality,
\be
E_Q \left[ \log q_S(s,t) \right] \geq \sup_y E_Q  E_y [ \int_0^{t-s} b(W(\tau),t-\tau)\,d\tau ] = 0.
\ee
Also, it follows from the Borell inequality and Theorem 3.2 of \cite{Adler2} (p. 63, let $\alpha = 1$) that there is a finite constant $K_0 > 0$ such that
\be
E_Q e^{\sup_y \int_0^{t-s} b(W^y(\tau),t-\tau)\,d\tau} < K_0 < \infty \label{e109}
\ee
if $\hat \sigma^2 < \frac{1}{2}$ . Thus by (\ref{e110}), there is a constant $K_1> 0$ such that (\ref{e109}) holds when $\abs{t - s} \leq K_1 $. Now for any $s < t$,  let $N$ be the smallest integer greater than $\abs{t-s}/K_1$ and $s = t_0 < t_1 < t_2 < \dots < t_N = t$ with $\abs{t_{i+1} - t_i} = \Delta t = \abs{t-s}/N \leq K_1$ for all $i = 0,\dots, N-1$. Using the subadditivity of $\log(q_S(s,t))$ and Jensen's inequality, we see that
\br
E_Q \log(q_S(s,t)) & \leq & \sum_{i=0}^{N-1} E_Q \log(q_S(t_i,t_{i+1})) \no \\
& \leq & \sum_{i=0}^{N-1} \log(E_P E_Q e^{\sup_y \int_0^{t_{i+1} - t_i} b(W^y(\tau),t_{i+1}-\tau)\,d\tau}). \no
\er
By (\ref{e109}), this right side is bounded by $N \log(K_0) < \infty$. Moreover,
\br
\sup_{s,t \in [0,K_1]} \abs{\log q_S(s,t)} \leq K_1 \sup_{t \in [0,K_1], y \in D} \abs{b(y,t)} \label{e114}
\er
which is integrable. Hence, we can apply the subadditive ergodic theorem to conclude that the limit
\be
\lim_{t \to \infty} \frac{\log(q_S(0,t))}{t} = \inf_t \frac{E_Q \left[\log( q_S(0,t))\right]}{t} = \rho_S, \label{e40}
\ee
holds almost surely with $\rho_S$ a constant, $\rho_S \in [0, \infty)$. Note that convergence along continuous time follows from (\ref{e114}), the continuity of $q_S(0,t)$ and Theorem 2.5 of \cite{AK}.

Clearly $\rho_I \leq \rho_S$. To show that $\rho_I = \rho_S$, we will need a kind of Harnack inequality to compare the quantities $q_I(0,t)$ and $q_S(0,t)$. Such a result has been obtained in \cite{CM} in the case that $b(y,t)$ is Gaussian in both space and time, with a white-noise temporal dependence.
Under the assumptions A5-A6, however,
the arguments of \cite{CM} imply that the following estimate also holds in the present case.
\begin{theo} \cite{CM} \label{theo:harnack}
For any fixed $M > 0$, there are positive constants $c_1, c_2$ such that outside an event of $Q-$probability $e^{-\frac{1}{4} n^{5/6}}$, one has
\be
\inf_{y \in D} E_y\left[ \xi(n,X) \right] \geq c_1 e^{-c_2 n^{11/12}} \left( \sup_{y \in D } E_y\left[ \xi(n,X) \right] - e^{-\frac{1}{4}n^{7/6}} \right). \no
\ee
\end{theo}

From this result it follows immediately that
\be
\lim_{n \to \infty} \frac{1}{n} \log \inf_{y} \phi(y,n) = \rho_I = \rho_S, \no
\ee
and by (\ref{e39}) and (\ref{e40}), we see that this extends to continuous time
\be
\lim_{t \to \infty} \frac{1}{t} \log \inf_{y} \phi(y,t) = \rho_I = \rho_S = \mu(\lambda) - f'(0) - \frac{\lambda^2}{2}. \label{e84}
\ee
We have now shown that $\rho(\lambda)$ (and thus $\mu(\lambda)$) in (\ref{e96}) is well-defined, independent of $z \in R^n$, almost surely with respect to $Q$.

The convexity of $\mu(\lambda)$ follows from the convexity of $\rho(\lambda)$. Let $r \in (0,1)$, $\lambda_1, \lambda_2 \in R$.  By H\"older's inequality,
\be
E[e^{r\lambda_1 I + (1-r)\lambda_2 I}] \leq E[e^{\lambda_1 I}]^r  E[e^{\lambda_2 I}]^{1-r}.\no
\ee
Applying this to (\ref{e96}), we conclude that
\be
\rho(r \lambda_1 + (1-r) \lambda_2) \leq r \rho(\lambda_1) + (1-r) \rho(\lambda_2). \no
\ee
Clearly $\rho(0) = 0$. Since $b(y,t) $ has the same distribution as $-b(y,t)$, we see that $\rho(-\lambda) = \rho(\lambda)$ for all $\lambda \in R$. Hence $\rho$ and $\mu$ are even functions of $\lambda$. Since $\rho(0) = 0$ and $\rho$ is convex and even, we conclude that
\be
\rho(\lambda) \geq 0\;\;\;\text{and}\;\;\; \mu(\lambda) \geq \lambda^2/2 + f'(0),\;\;\;\;\forall \lambda \in R.\label{e97}
\ee
This completes the proof of Lemma \ref{lem:muexist}.

We conclude this section by demonstrating that $\rho(\lambda)$ grows linearly with $\lambda$ as $\lambda \to \infty$.

\begin{prop}\label{prop:rholinear}
Let $b(y,t)$ have the form
\be
b(y,t) = \sum_{j=1}^{N} b^j_1(y) b^j_2(t) \label{e100}
\ee
where $b^j_1(y)$ are Lipschitz continuous and periodic in y, and $b^j_2(t)$ are stationary centered Gaussian fields such that the Assumptions A1-A6 are satisfied for $b(y,t)$. Then there is a constant $K > 0$ such that for $\lambda$ sufficiently large, $\rho(\lambda) \geq K \lambda$.
\end{prop}
{\bf Proof:}
In the case that $b(y,t)$ is a Gaussian field with white-noise time dependence, the authors of \cite{CM} studied the behavior of $\rho(\kappa)$ as $\kappa \to 0$, where $\kappa > 0$ is a diffusion constant (replace $\Delta/2$ with $\kappa \Delta/2$ in (\ref{e86})). Here we modify some of their ideas to treat the large advection limit when $b$ has the form (\ref{e100}).

Let $A_t$ be the set of continuous, piecewise linear functions $g$ such that $g(0) = 0$, $g$ is linear on the intervals $[it/k,(i+1)t/k]$, and
\be
g((i+1)t/k) - g(t/k) = \pm t/k , \;\;\;\;\; i = 0,1,2,\dots,k-1 \no
\ee
where $k >0$ is a large integer that will depend on $t$. Using the subadditivity arguments of \cite{CM} one can show that there is a constant $c_1 >0$ such that
\be
Q\left( \sup_{f \in A_t} \int_0^t b(f(s),t-s) \,ds \geq (c_1 - \epsilon) t \right) \geq 1 - \epsilon \label{e98}
\ee
if $k$ is chosen to be sufficiently large, depending on $t$. This implies that by choosing $\epsilon$ small there is a set of probability at least $1/2$ such that we can find $f(s) = f(s,\hat \omega) \in A_t$ satisfying
\be
\int_0^t b(f(s,\hat \omega),t-s, \hat \omega) \,ds \geq \frac{c_1 t}{2}, \label{e99}
\ee
and we expect that Brownian paths staying close to $f$ will make a significant contribution to the exponential in the definition of $\rho(\lambda)$. For a constant $\gamma > 0$ to be determined and $f \in A_t$, we let $B_t(f, \gamma)$ be the $\gamma$-neighborhood of $f$ in $C([0,t],R^{n-1})$:
\be
B_t(f,\gamma) = \left \{ X \in C([0,t],R^{n-1})\,|\;\;\norm{X - f}_{C^0} < \gamma \right \}. \no
\ee
Then there are constants $K_1, K_2$ independent of $t$ and of $f \in A_t$  such that
\be
P(B_t(f,\gamma)) \geq K_1 e^{-K_2 (1 + 1/\gamma^2) t} \no
\ee
for $t > 1$ (see \cite{CM}). Now using assumption A6, we see that for any path $X \in B_t(f,\gamma)$,
\be
\abs{\int_0^t b(X(s),t-s) \,ds - \int_0^t b(f(s),t-s) \,ds} < \gamma M  \sum_{j=1}^{N}  \int_0^t \abs{b^j_2(s)}\,ds \label{e101}
\ee
where $M$ is the maximum of the Lipschitz constants for the functions $\{ b^j_1(y) \}^N_{j=1}$. By (\ref{e99}) and (\ref{e101}) with $\epsilon>0$ sufficiently small, there is a set of $Q$-probability at least $1/2$ such that
\br
E_P\left[ e^{\lambda \int_0^t b(W^y_s(s),t-s) \,ds }\right] & \geq & E_P\left[ e^{\lambda \int_0^t b(W^y_s(s),t-s) \,ds } \chi_{B_t(f,\gamma)} \right] \no \\
& = &  e^{\lambda c_1 t/2} e^{ - \lambda \gamma M V} P(B_t(f,\gamma)) \no \\
& \geq & e^{\lambda c_1 t/2} e^{ - \lambda \gamma M V} K_1 e^{-K_2 (1 + 1/\gamma^2) t}, \label{e102}
\er
where $V =  \sum_{j=1}^{N}  \int_0^t \abs{b^j_2(s)}\,ds$ and $f \in A_t$ is chosen to satisfy (\ref{e99}).  For $t $ large, independently of $\lambda$ and $\gamma$, $V$ can be bounded by
\be
V \leq t (\sum_{j=1}^{N}  E[\abs{b^j_2(0}] + 1) \no
\ee
except on a set of probability less than $1/4$. Therefore, we can choose $\gamma$ small so that
\be
\gamma \leq \frac{c_1}{4 M (\sum_{j=1}^{N}  E[\abs{b^j_2(0}] + 1)}. \no
\ee
Then by choosing $\lambda$ large we obtain from (\ref{e102})
\be
E_P\left[ e^{\lambda \int_0^t b(W^y_s(s),t-s) \,ds }\right]  \geq  e^{\lambda c_1 t/8} \no
\ee
with $Q$-probability at least $1/4$, for $t$ sufficiently large, independently of $\lambda$. Since the limit defining $\rho(\lambda)$ exists $Q$-almost surely, this establishes the lemma with $K = c_1 /8$. \qed

\bibliographystyle{plain}

\begin{thebibliography}{99}
\bibitem{Adler2}R. Adler,
``An Introduction to Continuity, Extrema and Related Topics for General
Gaussian Processes", Institute of Math Stat, Lecture Notes-Monograph Series,
Vol. 12, 1990.

\bibitem{AK}
M.A. Akcoglu and U. Krengel,
{\em Ergodic theorems for superadditive processes},
J. Reine Angew Math. Vol. 323 (1981), pp. 53-67.

\bibitem{B2}H. Berestycki,
{\em The influence of advection on the propagation of fronts in reaction-diffusion equations},
in ``Nonlinear PDEs in Condensed Matter and Reactive Flows", NATO Science Series C, 569, H. Berestycki and Y. Pomeau
eds, Kluwer, Doordrecht, 2003.

\bibitem{BH1} H. Berestycki and F. Hamel,
{\em Front Propagation in Periodic Excitable Media},
Comm. in Pure and Appl. Math., Vol. 60, (2002), pp. 949-1032.

\bibitem{berestycki:tfc92}
H. Berestycki and L. Nirenberg, {\em Travelling fronts in cylinders},
  Ann. Inst. H. Poincar\'e Anal. Non Lin\'eaire, 9 (1992), pp.~497--572.

\bibitem{BHN}H. Berestycki, F. Hamel, N. Nadirashvili,
{\em Elliptic eigenvalue problems with large drift and applications to nonlinear propagation
phenomena}, Comm. Math Physics, 253(2), pp 451-480, 2005.

\bibitem{Min}R. Brent,``Module DFMIN in NMS",\hfill

http://gams.nist.gov/serve.cgi/Module/NMS/DFMIN/5671

\bibitem{Const1} P. Constantin, A. Kiselev, A. Oberman, L. Ryzhik,
{\em Bulk burning rate in passive-reactive diffusion}, Arch Rat.
Mech Analy, 154, (2000), pp. 53-91.

\bibitem{EL1}R.S. Ellis,
``Entropy, Large Deviations, and Statistical Mechanics",
Springer-Verlag: New York, 1985.

\bibitem{CarMol} R.A. Carmona and S.A. Molchanov,
{\em Parabolic Anderson problem and intermittency}. Mem. Amer. Math. Soc. 108 (1994), no. 518, viii+125

\bibitem{CD}J. Conlon, C. Doering,
{\em On Traveling Waves for the Stochastic FKKP Equation}, Jour. Stat Physics, to appear, 2005.

\bibitem{CM}M. Cranston and T. Mountford,
{\em Lyapunov Exponent for the Parabolic Anderson Model in $R^d$}, preprint, 2004.

\bibitem{CW} P. Clavin, and F. A. Williams,
{\em Theory of premixed-flame propagation in large-scale
turbulence}, J. Fluid Mech., 90, (1979), pp. 598-604.

%\bibitem{ES}
%{\sc L.C. Evans and P.E. Souganidis},
%{\em A PDE approach to geometric optics for certain semiliear parabolic equations},
%Indiana U. Math. J. 38, (1989), pp. 141-72.

\bibitem{FR1}M.I. Freidlin,
``Functional Integration and Partial Differential Equations", Ann. Math. Stud. 109,
Princeton University Press, Princeton, NJ, 1985.

\bibitem{FW1}M.I. Freidlin and A.D. Wentzell,
``Random Perturbations of Dynamical Systems". Springer-Verlag: New York, 1998.

\bibitem{HPS}S. Heinze, G. Papanicolaou, A. Stevens,
{\em Variational principles for propagation speeds in inhomogeneous media},
SIAM J. Applied Math, 62, no. 1, (2001), pp. 129 - 148.

\bibitem{Kato}T. Kato,
``Perturbation Theory for Linear Operators". Springer-Verlag: Berlin, 1995.

\bibitem{KR}A. Kiselev, L. Ryzhik,
{\em Enhancement of the traveling front speeds in reaction-diffusion
equations with advection}, Ann. de l'Inst. Henri Poincar\'e,
Analyse Nonlin\'eaire, 18, (2001), pp. 309--358.

\bibitem{King}J. P. C. Kingman,
{\em The Ergodic Theory of Subadditive Stochastic Processes}, Journal of the Royal Statistical Society, Series B, Vol. 30, No. 3, (1968), pp. 499-510.

\bibitem{Kh}B. Khouider, A. Bourlioux, A. Majda,
{\em Parameterizing turbulent flame speed-Part I: unsteady shears, flame
residence time and bending}, Combustion Theory and Modeling, 5 (2001),
pp. 295-318.

\bibitem{KoPl}P. Kloeden and E. Platen,
``Numerical Solution of Stochastic Differential Equations",
Springer-Verlag, 1999.

%\bibitem{KZ1}
%{\sc A. Kiselev and A. Zlatos}
%{\em Quenching of combustion by shear flows}, preprint. To appear in Duke Math. J., 2005.

%\bibitem{PLS}{\sc P.-L. Lions and P.E. Souganidis}
%{\em Homogenization of Viscous Hamilton-Jacobi Equations in Stationary Ergodic Media}, Comm. Partial Diff. Eqn. 30, (2005), pp. 335-375.

\bibitem{MS1}A. Majda and P.E. Souganidis,
{\em Large scale front dynamics for turbulent reaction-diffusion equations
with separated velocity scales}, Nonlinearity, 7 (1994), pp. 1-30.

\bibitem{MS2}A. Majda and P.E. Souganidis,
{\em Flame fronts in a turbulent combustion model with fractal velocity
fields}, Comm Pure Appl Math, Vol. LI (1998), pp. 1337-1348.

\bibitem{MShen}J. Mierczynski, W. Shen,
{\em Exponential separation and principle Lyapunov exponent/spectrum for
random/nonautonomous parabolic equations}, J. Differential Equations, 191(2003), pp 175-205.

\bibitem{MS}C. Mueller, R. Sowers,
{\em Random Traveling Waves for the KPP equation with Noise},
J. Functional Analysis, 128(1995), pp  439-498.

\bibitem{NRX}J. Nolen, M. Rudd, J. Xin,
{\em Existence of KPP fronts in spatially-temporally periodic advection and
variational principle for propagation speeds}, Dynamics of PDE, Vol. 2, No. 1, pp 1-24, 2005.

\bibitem{NX1}J. Nolen, J. Xin,
{\em Reaction diffusion front speeds in spatially-temporally periodic
shear flows}, SIAM J. Multiscale Modeling and Simulation,
Vol. 1, (2003), No. 4, pp. 554-570.

\bibitem{nolen:ekt04}
J. Nolen, J. Xin, {\em Existence of {KPP} type fronts in space--time
  periodic shear flows and a study of minimal speeds based on variational
  principle}, Discrete and Continuous Dynamical Systems, Vol. 13, No. 5 (2005), pp. 1217-1234.

\bibitem{NX2}J. Nolen, J. Xin,
{\em A Variational Principle Based Study of KPP Minimal
Front Speeds in Random Shears}, Nonlinearity 18 (2005), pp 1655 - 1675.

\bibitem{NX3}J. Nolen, J. Xin,
{\em Min-Max Variational Principle and Front Speeds in Random Shear Flows},
math.AP/0501445, www.arxiv.org, 2005, to appear in Methods and Applications of Analysis.

\bibitem{Peters}N. Peters,
{\em Turbulent Combustion}, Cambridge University Press, 2000.

\bibitem{Ro}P. Ronney,
{\em Some open issues in premixed turbulent combustion},
in: Modeling in Combustion Science (J. D. Buckmaster and T. Takeno, Eds.),
Lecture Notes In Physics, Vol. 449, Springer-Verlag, Berlin, (1995), pp. 3-22.

\bibitem{Shen}W. Shen,
{\em Traveling Waves in Diffusive Random Media}, J. Dynamics Diff. Eqs. 16(2004), No. 4,
pp 1011-1060.

\bibitem{Vlad}N. Vladimirova, P. Constantin, A. Kiselev, O. Ruchayskiy, L. Ryzhik,
{\em Flame enhancement and quenching in fluid flows}, Combust. Theory and Modeling,
7, (2003), pp. 487-508.

\bibitem{xin:est91}
J. Xin, {\em Existence and stability of travelling waves in periodic
  media governed by a bistable nonlinearity}, J. Dynamics Diff. Eqs., 3 (1991),
  pp.~541--573.

\bibitem{xin:epf92}
J. Xin, {\em Existence of planar
  flame fronts in convective--diffusive periodic media}, Arch. Rat. Mech.
  Anal., 121 (1992), pp.~205--233.

\bibitem{Xin1}
J. Xin,
{\em Front propagation in heterogeneous media,}
SIAM Review, Vol. 42, No. 2, June 2000, pp. 161-230.

\bibitem{Xin3}J. Xin,
{\em KPP front speeds in random shears and the parabolic Anderson problem},
Methods and Applications of Analysis, Vol. 10, No. 2, (2003), pp. 191-198.

\bibitem{Yak}V. Yakhot,
{\em Propagation velocity of premixed turbulent flames},
Comb. Sci. Tech 60, (1988), p. 191.

%\bibitem{ZL1}
%{\sc A. Zlatos},
%{\em Quenching and propagation of combustion without ignition temperature cutoff},
%Nonlinearity, Vol. 18, (2005), pp. 1463-1475.


\end{thebibliography}

\end{document}